\documentclass[a4paper]{scrartcl}
\usepackage[utf8]{inputenc}
\usepackage{amsmath,amssymb,amsthm,amsfonts}
\usepackage{mathrsfs}
\usepackage{mathtools,thmtools}
\usepackage[british]{babel}
\usepackage{csquotes}
\usepackage{hyperref}
\usepackage[nameinlink,noabbrev,capitalise]{cleveref}
\usepackage{enumerate}
\usepackage{xcolor}
\usepackage{comment}
\usepackage{romanbar}
\usepackage{bm}

\RequirePackage{tikz}
    \usetikzlibrary{positioning}
    \usetikzlibrary{decorations.pathreplacing}

\usepackage[maxbibnames=5,sorting=nyt,sortcites,giveninits=true]{biblatex}
\DeclareNameAlias{sortname}{family-given}
\usepackage[T2A,T1]{fontenc}
\DeclareSymbolFont{cyrillic}{T2A}{cmr}{m}{n} 
\DeclareMathSymbol{\rusp}{\mathalpha}{cyrillic}{239}
\DeclareMathSymbol{\rusz}{\mathalpha}{cyrillic}{214}
\DeclareMathSymbol{\rusg}{\mathalpha}{cyrillic}{227}
\DeclareMathSymbol{\rush}{\mathalpha}{cyrillic}{237}
\DeclareMathSymbol{\rusd}{\mathalpha}{cyrillic}{196}

\newcommand{\R}     {\mathbb{R}}
\newcommand{\N}     {\mathcal{N}}

\renewcommand{\P}   {\mathbb{P}}
\newcommand{\E}     {\mathbb{E}}

\newcommand{\eps}   {\varepsilon}
\let\S\relax
\newcommand{\S}     {\mathcal{S}}
\newcommand{\calP}     {\mathcal{P}}

\newcommand{\TV}    {\textup{TV}}
\newcommand{\KL}    {\textup{KL}}

\newcommand{\tmix}{t_{\operatorname{mix}}}
\newcommand{\ind}{\mathbf{1}}

\newcommand{\comp}{\mathsf{c}}

\DeclarePairedDelimiterX{\norm}[1]{\lVert}{\rVert}{#1}

\DeclareMathOperator{\Unif} {Unif}

\DeclareMathOperator{\Law}  {Law}

\theoremstyle{plain}
\newtheorem{theorem}{Theorem}
\newtheorem{lemma}[theorem]{Lemma}
\newtheorem{proposition}[theorem]{Proposition}

\theoremstyle{definition}

\newtheorem{definition}{Definition}

\newtheorem{remark}{Remark}
\newtheorem{assumption}{Assumption}

\title{High-probability Uniform Mixing via Dobrushin Contraction and Accelerated Hamiltonian Monte Carlo}
\author{%
  \Large Stefan Oberd\"orster%
  \thanks{%
    E-mail:\ \href{mailto:oberdoerster@uni-bonn.de}
                 {\texttt{oberdoerster@uni-bonn.de}},
    ORCID:\ \href{https://orcid.org/0009-0004-0734-0101}
                {0009-0004-0734-0101}.%
  }\\[0.6em]
  \large Institute for Applied Mathematics, University of Bonn
}
\date{\normalsize\today}

\begin{document}

\maketitle

\begin{abstract}
\noindent
Hamiltonian Monte Carlo, since its invention and subsequent rise to widespread use, has been motivated by its ability to accelerate convergence to equilibrium compared to diffusive methods, both in terms of dependence on condition number and dimension.
In the present work, we establish this acceleration simultaneously in both aspects for critically damped Metropolis-adjusted kinetic Langevin with persistent uniforms.
More precisely, for target measures $\mu(dx)\propto e^{-U(x)}dx$ with $m$-strongly convex, $L$-gradient Lipschitz, and $L_H$-Frobenius Hessian Lipschitz potentials, we show that the method mixes to $\eps$-accuracy in total variation using $\mathcal O\bigl((\kappa^{1/2}+L_H^{1/2}m^{-3/4})\,d^{1/4}\bigr)$ gradient evaluations up to logarithmic factors, both uniformly over starting points in domains of high equilibrium probability and from warm starts.
The proof incorporates ideas from previous studies \cite{Nawaf26persistence,chen2023} that show either aspect of acceleration separately, into novel coupling technology.
This new coupling technology for general Markov chains extends earlier frameworks aimed at capturing fast local equilibration in the presence of slow global mixing.
Specifically, via a weak-$(1,1)$ maximal inequality that follows from Hopf's maximal ergodic theorem, we uncover a stability mechanism in high-probability domains with respect to equilibrium.
This enables straightforward localization to naturally arising domains of fast local mixing, and seamlessly fits broad local frameworks via Dobrushin contraction \cite{DPMC}, Wasserstein contraction \cite{Eberle14} and local Harris' theorems \cite{Yang23}.
We develop the Dobrushin formulation.
\end{abstract}

\section{Introduction}

Consider the Markov chain Monte Carlo problem of sampling a target probability distribution
\[ \mu(dx)\ \propto\ e^{-U(x)}\,dx\quad\text{on $\R^d$} \]
with potential $U:\R^d\to\R$.
Hamiltonian Monte Carlo methods, since their invention in lattice field theory \cite{Duane87,Horowitz91}, aim at improving performance over previously used diffusive methods \cite{Rossky78,Pangali78}.
Hamiltonian Monte Carlo discretizes the kinetic Langevin diffusion
\begin{equation}\label{eq:LD}
    dX_t\ =\ V_t\,dt,\quad dV_t\ =\ -\nabla U(X_t)\,dt -\gamma V_t\,dt +\sqrt{2\gamma}\,dB_t,
\end{equation}
where $(B_t)_{t\geq0}$ is a standard Brownian motion on $\R^d$ and $\gamma>0$ quantifies friction.
This is achieved through alternating between Hamiltonian flow with respect to 
\[ H(x,v)\ =\ U(x)\ +\ \tfrac12|v|^2 \]
and refreshing the velocity from $\mathcal N(0,I_d)$, either partially at rate $\gamma$ or fully after integration times of order $\gamma^{-1}$.
For implementability, Hamiltonian flow is approximated by a numerical scheme, most commonly the Størmer--Verlet, aka. leapfrog, integrator \cite{Verlet67,Strang68}.
A Metropolis--Hastings correction \cite{Metropolis53,Hastings70,Tierney98} eliminates discretization bias resulting from integration and ensures the methods to leave invariant the lifted target
\begin{equation}\label{eq:hatmu}
    \hat\mu(dx\,dv)\ =\ \bigl(\mu\otimes\mathcal N(0,I_d)\bigr)(dx\,dv)\ \propto\ e^{-H(x,v)}\,dx\,dv\quad\text{on phase space $\R^{2d}$.}
\end{equation}

The superior practical performance of Hamiltonian Monte Carlo, as opposed to diffusive methods, is rooted in the persistence of momentum on time scales of order $\gamma^{-1}$.
Endowing the dynamics with the ability of directed motion, this has two important implications.
For their intuition, assume $\mu$ to be a product measure whose univariate components' scales range from $L^{-1/2}$ to $m^{-1/2}$ for $0<m\leq L$.
The most simple example is given by the Gaussian $\mathcal N(0,A^{-1})$ with $m\,I_d\leq A\leq L\,I_d$.
\begin{enumerate}[(i)]
\item
On the one hand, at critical damping $\gamma^{-1}\asymp m^{-1/2}$, the directed motion persists over time comparable to the target's largest scale, enabling ballistic exploration in time $t_{\mathrm{bal}}\asymp m^{-1/2}$.
In contrast, diffusive motion covers distance $t^{1/2}$ in time $t$ and hence requires time $t_{\mathrm{dif}}\asymp m^{-1}$ to explore the target distribution's largest scale.
\item
On the other hand, heuristically, as discussed in \cite[Sec.~6]{Nawaf18}, the equilibrium mean discretization error of the leapfrog integrator with step size $h>0$ is bounded by $\mathcal O(h^4L^2d)$, independently of the integration time, and hence controlled for $h^{-1}\asymp L^{1/2}d^{1/4}$.
Along the critical ballistic motion, this results in $t_{\mathrm{bal}}/h\asymp \kappa^{1/2}d^{1/4}$ leapfrog steps, where $\kappa=L/m$ denotes the condition number.
In comparison, the Metropolis-adjusted Langevin algorithm, which makes one leapfrog step per iteration approximating an $\asymp h^2$ time step of the overdamped Langevin diffusion, explores space diffusively, resulting in $t_{\mathrm{dif}}/h^2\asymp\kappa\,d^{1/2}$ steps.
\end{enumerate}
\noindent
The suggested improvement from $\kappa\,d^{1/2}$ to $\kappa^{1/2}d^{1/4}$ leapfrog steps is known as \emph{diffusive-to-ballistic acceleration}.
It motivated Hamiltonian Monte Carlo since its invention \cite{Kennedy96,Gupta88}.
Note that acceleration impacts both the convergence of the ideal dynamics, such as \eqref{eq:LD}, in (i), as well as discretization in (ii).

On the level of ideal dynamics, accelerated convergence to equilibrium under general assumptions was first proved in $L^2$ in \cite{CaoLuWang23,lu2022explicit} via the variational hypocoercivity framework based on space-time Poincar\'e inequalities of \cite{albritton2024variational}.
Subsequently, \cite{EberleLoerler24} cast the analysis in the language of lifts and showed the square-root acceleration to be the best achievable through lifting.
Following developments in $L^2$, \cite{Lu26entropy} extended the hypocoercivity approach to entropy while retaining acceleration, see also \cite{Wibisono26}.

All of the aforementioned works are concerned with ideal dynamics and neither cover discretization nor Metropolis-adjustment.
Acceleration in terms of the explicit dimension scaling under general assumptions was first established in \cite{chen2023} via conductance techniques.
Being limited to short persistent motion, on time scales $L^{-1/2}$ in the above setting, inducing diffusive exploration in the sense of (i), their study does not capture acceleration in the condition number.
Extending the entropic hypocoercivity for ideal dynamics of \cite{Lu26entropy} to a discretization of \eqref{eq:LD}, \cite{Nawaf26persistence} showed condition number acceleration for a discretized and Metropolis-adjusted method for the first time under general assumption, reporting a non-accelerated dimension dependence however.

Notably, no previous work simultaneously captures both acceleration in condition number and dimension, jointly motivating Hamiltonian Monte Carlo, under general assumptions.
Our first main result \Cref{thm:main} fills this gap by incorporating ideas from \cite{Nawaf26persistence} and \cite{chen2023} into novel general-purpose coupling technology---our second main result discussed below.
More precisely, assuming the potential $U$ to be $m$-strongly convex, $L$-gradient Lipschitz, and $L_H$-Hessian-Lipschitz with respect to the Frobenius norm, we establish that the Metropolis-adjusted OBABO discretization of \eqref{eq:LD} using persistent uniforms \cite{Neal20} at critical friction $\gamma=2\sqrt m$ mixes to accuracy $\eps\in(0,1/2)$ in total variation in
\[ \widetilde{\mathcal O}\left(\bigl(\kappa^{1/2}+L_H^{1/2}m^{-3/4}\bigr)\bigl(d+\log(\delta\eps)^{-1}\bigr)^{1/4}\log^{1/2}(\delta\eps)^{-1}\log\eps^{-1}\right) \]
steps, for any $\delta\in(0,1/2)$, both
\begin{enumerate}[(i)]
\item uniformly over starting points in suitable high-probability domains $D\subseteq\R^{2d}$ such that $\hat\mu(D)\geq1-\delta$, and
\item uniformly from $L^p$-warm starts $\hat\nu\ll\hat\mu$ such that $\beta_p(\hat\nu):=\bigl\|\frac{d\hat\nu}{d\hat\mu}\bigr\|_{L^p(\hat\mu)}\leq\delta^{-1}$ for fixed $p\in(1,\infty]$,
\end{enumerate}
where $\widetilde{\mathcal O}$ omits logarithmic factors in the model parameters and double-logarithmic ones in $(\eps\delta)^{-1}$.

In particular, for $\log(\eps\delta)^{-1}=\widetilde{\mathcal O}(1)$, allowing the inverse accuracy, inverse exceptional mass, and initial warmness to depend polynomially on the model parameters, we show mixing in
\[ \widetilde{\mathcal O}\left(\bigl(\kappa^{1/2}+L_H^{1/2}m^{-3/4}\bigr)\,d^{1/4}\right) \]
steps.
In this setting, we improve upon \cite{chen2023} and \cite{Nawaf26persistence}, showing mixing with complexity $\widetilde{\mathcal O}\bigl((\kappa+L_H^{2/3}m^{-1})\,d^{1/4}\bigr)$ from $L^\infty$-warm starts and $\widetilde{\mathcal O}\bigl((\kappa d)^{1/2}+(L_H^{1/2}m^{-3/4})\,d^{3/4}\bigr)$ from $L^2$-warm starts, respectively.
It is important to emphasize however, that both previous works consider weaker assumptions, \cite{chen2023} replacing strong log-concavity of $\mu$ by an isoperimetric inequality, and \cite{Nawaf26persistence} the Frobenius norm in the Hessian-Lipschitzness of $U$ by the spectral norm.

Moreover, our results permit deterministic initialization in suitable domains of large equilibrium mass, which extends the scope regarding initialization beyond previous studies.
Specifically, critically tuned Metropolis-adjusted Markov chains, among others, exhibit highly local mixing, featuring coexistence of fast local mixing in typical equilibrium states with slow global mixing in the tails \cite{Diss}.
Previous analyses capturing fast local mixing assume warm initialization \cite{chen22,chen2023,Chen23,Chewi21}, while previous studies allowing non-warm initialization are limited to slow global mixing \cite{Eberle14,Eberle20,BouRabeeOberdoerster2024,Nawaf13,Alain22} or special cases \cite{Diss,BouRabeeOberdoersterNUTS1,NUTS2}.
Even though several localized approaches accommodating non-warm initialization are available, e.g., via Dobrushin contraction \cite{DPMC}, Wasserstein contraction \cite{Eberle14}, or Harris' theorems \cite{Yang23}, an easily applicable stability mechanism was lacking.

We contribute such a stability mechanism in high-probability domains with respect to the invariant measure.
It is based on a weak-$(1,1)$ maximal inequality that follows from Hopf's maximal ergodic theorem.
Our second main result \Cref{thm:mixlocDobHopf} combines this stability mechanism with the localized Dobrushin contraction framework of \cite{DPMC}, creating a general-purpose approach to \emph{uniform geometric mixing in high-probability domains}.
It enables fast local mixing guarantees under mild assumptions on initialization.
For instance, it yields the first study of Hamiltonian Monte Carlo capturing fast local mixing under general assumptions from deterministic starts.
More broadly, the stability mechanism immediately applies to the other localized approaches \cite{Eberle14,Yang23} and their variants \cite{BouRabeeOberdoerster2024,Diss,Nawaf26persistence}.

\paragraph{Structure of the Paper}
We state our main results \Cref{thm:main,thm:mixlocDobHopf}, on the mixing of Metropolis-adjusted kinetic Langevin and the underlying general framework for uniform mixing in high-probability domains, in \Cref{sec:main}.
Then, we first prove the general framework in \Cref{sec:locDobmix}, followed by the details of its application to Metropolis-adjusted kinetic Langevin in \Cref{sec:proofs}.

\section{Main Results}\label{sec:main}

We begin by introducing Metropolis-adjusted kinetic Langevin with persistent uniforms and some general notation in \Cref{sec:MAKLpu,sec:not}, followed by the statement of \Cref{thm:main} in \Cref{sec:MAKLpuMix} and \Cref{thm:mixlocDobHopf} in \Cref{sec:mixlocDobHopf}.

\subsection{Metropolis-adjusted Kinetic Langevin with Persistent Uniforms}\label{sec:MAKLpu}

Metropolis-adjusted kinetic Langevin is an instance of Hamiltonian Monte Carlo employing the OBABO scheme \cite{Bussi07,Leimkuhler13} to discretize the kinetic Langevin diffusion \eqref{eq:LD} with friction $\gamma>0$.
One step of size $h>0$ of OBABO takes the palindromic form 
\[ \Phi_{h,\gamma}^{\xi,\xi'}\ =\ O_{\gamma,h/2}^{\xi'}\circ\theta_{h/2}^{(B)}\circ\theta_h^{(A)}\circ\theta_{h/2}^{(B)}\circ O_{\gamma,h/2}^{\xi} \]
involving discrete kinetic and potential flows
\[ \theta_h^{(A)}(x,v)\ =\ (x+hv,\,v)\quad\text{and}\quad\theta_h^{(B)}(x,v) \ = \ \bigl(x,\,v-h\nabla U(x)\bigr), \]
as well as two Ornstein--Uhlenbeck (OU) half-steps of
\[ O_{\gamma,h}^{\zeta}(x,v) \ = \ \bigl(x,\,e^{-\gamma h}v+(1-e^{-2\gamma h})^{1/2}\zeta\bigr) \]
partially refreshing velocity with independent $\xi,\xi'\sim\mathcal N(0,I_d)$.
The kinetic and potential flow portion between partial refreshes corresponds one step of the leapfrog integrator \cite{Verlet67,Strang68}, which we denote by $\Phi_h$.

As velocity refreshment from $\mathcal N(0,I_d)$ preserves the lifted target on phase space $\hat\mu$, see \eqref{eq:hatmu}, only the leapfrog step incurs discretization error, measured in terms of the energy error $\Delta H_h=H\circ\Phi_h-H$.
The resulting bias is corrected via a Metropolis--Hastings mechanism \cite{Metropolis53,Hastings70}.
More precisely, the deterministic proposal on phase space $S\circ\Phi_h$ with $S(x,v)=(x,-v)$, due to being a volume-preserving involution \cite{Nawaf18}, with acceptance probability $e^{-\Delta H_h^+}$ yields a well-defined Metropolis--Hastings transition leaving $\hat\mu$ invariant \cite{Tierney98}.
To keep persistent momentum across Metropolis-adjusted OBABO steps, an additional velocity flip is added after the Metropolis--Hastings accept--reject step.

One complete transition from $Z_{k-1}$ to $Z_k$ of \emph{Metropolis-adjusted kinetic Langevin} takes the form
\[ Z_k\ =\ \begin{cases}
    \Phi_{h,\gamma}^{\xi_k,\xi'_k}(Z_{k-1}) & \text{if $U<e^{-\Delta H_h\circ O_{\gamma,h/2}^{\xi_k}(Z_{k-1})}$,} \\[1ex]
    O^{\xi'_k}_{\gamma,h/2}\circ S\circ O^{\xi_k}_{\gamma,h/2}(Z_{k-1}) & \text{otherwise,} 
    \end{cases} \]
where $\xi_k,\xi'_k\sim\mathcal N(0,I_d)$ and $U\sim\Unif(0,1)$ are independent and independent of $Z_{k-1}$.
Note that only the leapfrog step is Metropolis-adjusted, resulting in the acceptance probability being evaluated after the first OU half-step, and that the additional velocity flip remains in case of rejection.

As motivated in the introduction, the desirable properties of Hamiltonian Monte Carlo arise through the persistence of momentum.
However, upon rejection, the velocity is flipped, breaking persistent motion.
Therefore, it is desirable for the occurrence of rejection to be controlled over \emph{epochs} of order $(\gamma h)^{-1}$ transitions of Metropolis-adjusted kinetic Langevin, corresponding to physical time scales of order $\gamma^{-1}$.
When independent uniforms are used across transitions, the probability of rejection accumulates linearly over the epoch.
In contrast, Hamiltonian Monte Carlo methods, such as Randomized Hamiltonian Monte Carlo \cite{BoSa2017}, refreshing the velocity fully after integrating and Metropolis-adjusting Hamiltonian flow over $(\gamma h)^{-1}$ leapfrog steps, do not suffer this accumulation in the same way.
Instead, part of the leapfrog energy error telescopes over several consecutive steps, resulting in cancellation.

To supply Metropolis-adjusted kinetic Langevin with similar cancellation and therewith better behaved rejection probabilities over epochs of transitions, \cite{Neal20} introduced non-reversibly updated persistent uniforms.
Specifically, over an epoch of transitions, instead of refreshing the uniform random variables in each transition, they are deterministically updated upon acceptance and kept upon rejection.
This results in the following variant of Metropolis-adjusted kinetic Langevin, studied in \cite{Nawaf26persistence} and this work.

\begin{definition}[Metropolis-adjusted Kinetic Langevin with Persistent Uniforms]\label{def:MAKLApu}
An epoch of $E\in\mathbb N$ transitions of \emph{Metropolis-adjusted kinetic Langevin with persistent uniforms} with friction $\gamma>0$ and step size $h>0$ from $z=Z_0\in\R^{2d}$ proceeds as follows:
\begin{enumerate}[(i)]
\item Draw independent $U_0\sim\Unif(0,1)$ and $\xi_1,\xi'_1,\dots,\xi_E,\xi'_E\sim\N(0,I_d)$,
\item For $k=1,\dots,E$, iteratively set
\end{enumerate}
\[ (Z_k,U_k) = \left\{\begin{array}{@{}ll@{}}
    \Bigl(\Phi_{h,\gamma}^{\xi_k,\xi'_k}(Z_{k-1}),\ U_{k-1}e^{\Delta H_h\circ O_{\gamma,h/2}^{\xi_k}(Z_{k-1})}\Bigr) & \text{if $U_{k-1}<e^{-\Delta H_h\circ O_{\gamma,h/2}^{\xi_k}(Z_{k-1})}$,} \\[1.5ex]
    \Bigl(O^{\xi'_k}_{\gamma,h/2}\circ S\circ O^{\xi_k}_{\gamma,h/2}(Z_{k-1}),\ U_{k-1}\Bigr) & \text{otherwise.} 
    \end{array}\right. \]
We denote the Markov transition kernel describing such an epoch of steps by
\[ \Pi(z,\cdot)\ =\ \Law(Z_E) \]
and the acceptance events in the individual steps by
\[ A_k\ =\ \Bigl\{U_{k-1}<e^{-\Delta H_h\circ O_{\gamma,h/2}^{\xi_k}(Z_{k-1})}\Bigr\}\quad\text{for $1\leq k\leq E$.} \]
In case all steps are accepted, i.e., in the event $\bigcap_{k=1}^EA_k$, the epoch of Metropolis-adjusted steps reduces to the underlying epoch of OBABO steps, whose transition kernel we denote by
\[ \Pi_u(z,\cdot)\ =\ \Law\left(\Phi_{h,\gamma}^{\xi_E,\xi'_E}\circ\cdots\circ\Phi_{h,\gamma}^{\xi_1,\xi'_1}(z)\right). \]
\end{definition}

As shown in \cite[Lem.~2.3]{Nawaf26persistence}, $\Pi$ leaves $\hat\mu$ invariant.
As intended, the persistent uniforms provide the epoch of transitions with cancellation due to approximate leapfrog error telescoping.
More precisely, on the event $\bigcap_{k=1}^{n-1}A_k$ that all transitions up to the $n$-th were accepted, the subsequent acceptance event is
\[ A_n\ =\ \Bigl\{U_0<e^{-\sum_{k=1}^n\Delta H_h\circ O_{\gamma,h/2}^{\xi_k}(Z_{k-1})}\Bigr\}, \]
i.e., the leapfrog energy errors telescope up to OU half-steps.

\subsection{Notation}\label{sec:not}

Let us fix some notation that is used throughout the presentation.
Our general results are formulated for Markov chains with transition kernels $\pi$ and invariant probability measures $\mu$ on Polish state spaces $\S$, endowed with their $\sigma$-algebra $\mathfrak B$.
We denote the set of all probability measures on $(\S,\mathfrak B)$ by $\mathcal P(\S)$.
The results for Metropolis-adjusted kinetic Langevin specialize this setting to Euclidean space.

We quantify convergence to equilibrium in terms of the total variation distance, defined between any $\eta,\tilde\eta\in\mathcal P(\S)$ by $\TV(\eta,\tilde\eta)=\sup_{B\in\mathfrak B}|\eta(B)-\tilde\eta(B)|$.
Correspondingly, we define the mixing time to accuracy $\eps>0$ from start in $\nu\in\mathcal P(\S)$ by
\[ \tmix(\eps,\nu)\ =\ \inf\bigl\{n\geq0:\TV(\nu\pi^n,\mu)\leq\eps\bigr\}. \]
For any domain $D\subseteq\S$, we further define the uniform mixing time from start in $D$ by
\[ \tmix(\eps,D)\ =\ \sup\nolimits_{x\in D}\tmix(\eps,\delta_x). \]
For $p\in(1,\infty]$, the $L^p$-warmness of an initial distribution $\nu\in\mathcal P(\S)$ is defined by $\beta_p(\nu)=\bigl\|\frac{d\nu}{d\mu}\bigr\|_{L^p(\mu)}$ if $\nu\ll\mu$ and to be infinite otherwise.
We say $\nu$ to be $L^p$-warm if $\beta_p(\nu)<\infty$.

We write $a\lesssim b$ or $a=\mathcal O(b)$ if there exists an absolute constant $C>0$ such that $a\leq Cb$, and $a\asymp b$ if both $a\lesssim b$ and $b\lesssim a$.

\subsection{Mixing of Metropolis-adjusted Kinetic Langevin}\label{sec:MAKLpuMix}

We make the following assumption on the target distribution.

\begin{assumption}\label{ass:1}
The probability distribution $\mu$ on $\R^d$ admits a Lebesgue density proportional to $\exp(-U(x))$ with non-negative potential $U\in C^2(\R^d)$.
Further:
\begin{enumerate}[(i)]
\item
\emph{Strong Convexity and Gradient Lipschitzness:}
There exist $0<m\leq L$ such that
\[ m\,I_d\ \leq\ \nabla^2 U(x)\ \leq\ L\,I_d\quad\text{for all $x\in\R^d$.} \]
Under strong convexity, $U$ has a unique minimum, which we assume to be $U(0)=0$.
We denote the condition number by $\kappa=L/m$.
\item
\emph{Frobenius Hessian Lipschitzness:}
There exists $L_H\geq0$ such that
\[ \bigl\|\nabla^2U(x)-\nabla^2U(y)\bigr\|_F\ \leq\ L_H|x-y|\quad\text{for all $x,y\in\R^d$.} \]
\end{enumerate}
\end{assumption}

Frobenius Hessian Lipschitzness strengthens the also common assumption of Hessian Lipschitzness with respect to the spectral norm.
Note that if $U$ is $L_{H,s}$-Hessian Lipschitz with respect to the spectral norm, \Cref{ass:1}(ii) holds with $L_H=d^{1/2}L_{H,s}$.

\begin{theorem}[High-probability Uniform Mixing of Metropolis-adjusted Kinetic Langevin]\label{thm:main}
Suppose $\mu$ satisfies \Cref{ass:1} and, as in \Cref{def:MAKLApu}, let $\Pi$ be the transition kernel of an epoch of $E\in\mathbb N$ steps of Metropolis-adjusted kinetic Langevin with persistent uniforms with friction $\gamma=2\sqrt m$ and step size $h>0$.
Let $\eps,\delta\in(0,1/2)$.
Then, there exist a measurable $D\subseteq\R^{2d}$ and suitable $E,h$ such that
\begin{equation}\label{eq:mixtimeEpoch}
    \tmix(\eps,D)\ =\ \mathcal O(\log\eps^{-1})\quad\text{and}\quad\hat\mu(D^\comp)\ \leq\ \delta.
\end{equation}
The resulting complexity in terms of potential gradient evaluations to achieve $\eps$-accuracy in total variation uniformly from start in $D$ is
\begin{equation}\label{eq:cc}
    \mathcal O\left(\Lambda^{9/4}\bigl(\kappa^{1/2}+L_H^{1/2}m^{-3/4}\bigr)\bigl(d+\log(\delta\eps)^{-1}\bigr)^{1/4}\log^{1/2}(\delta\eps)^{-1}\log\eps^{-1}\right),
\end{equation}
where $\Lambda=\log\max\bigl(e,\kappa,L_H,L_H/m^{3/2},d,\log(\delta\eps)^{-1}\bigr)$.
\end{theorem}

The theorem guarantees uniform mixing from start in a domain of $\hat\mu$-mass at least $1-\delta$, with gradient complexity growing only logarithmically in $\delta^{-1}$ as $\delta\to0$.

\begin{remark}[Mixing from Warm Start]\label{rem:mixwarm}
Beyond uniformly in high-probability domains, we additionally show mixing of Metropolis-adjusted kinetic Langevin with persistent uniforms from warm starts under the same \Cref{ass:1} and friction $\gamma=2\sqrt m$.
From start in any $L^p$-warm initial distribution $\hat\nu\in\mathcal P(\R^{2d})$ for $p\in(1,\infty]$, \Cref{thm:mainwarm} below guarantees mixing to accuracy $\eps\in(0,1/2)$ for suitable $h,E$.
The resulting complexity in terms of potential gradient evaluations is
\[ \mathcal O\left(\Lambda^{9/4}\bigl(\kappa^{1/2}+L_H^{1/2}m^{-3/4}\bigr)\Bigl(d+\tfrac{p}{p-1}\log\bigl(\beta_p(\hat\nu)/\eps\bigr)\Bigr)^{1/4}\Bigl(\tfrac{p}{p-1}\log\bigl(\beta_p(\hat\nu)/\eps\bigr)\Bigr)^{1/2}\log\eps^{-1}\right), \]
where $\Lambda=\log\max\bigl(e,\kappa,L_H,L_H/m^{3/2},d,\frac{p}{p-1}\log\bigl(\beta_p(\hat\nu)/\eps\bigr)\bigr)$.
\end{remark}

As stated in the introduction, under \Cref{ass:1}, \Cref{thm:main} and \Cref{rem:mixwarm} together guarantee critically tuned Metropolis-adjusted kinetic Langevin with persistent uniforms to achieve accuracy $\eps\in(0,1/2)$ in total variation using
\[ \widetilde{\mathcal O}\left(\bigl(\kappa^{1/2}+L_H^{1/2}m^{-3/4}\bigr)\bigl(d+\log(\delta\eps)^{-1}\bigr)^{1/4}\log^{1/2}(\delta\eps)^{-1}\log\eps^{-1}\right) \]
potential gradient evaluations, for any $\delta\in(0,1/2)$, both
\begin{enumerate}[(i)]
\item uniformly over starting points in suitable high-probability domains $D\subseteq\R^{2d}$ such that $\hat\mu(D)\geq1-\delta$, and
\item uniformly from $L^p$-warm starts $\hat\nu\ll\hat\mu$ such that $\beta_p(\hat\nu):=\bigl\|\frac{d\hat\nu}{d\hat\mu}\bigr\|_{L^p(\hat\mu)}\leq\delta^{-1}$ for fixed $p\in(1,\infty]$,
\end{enumerate}
where $\widetilde{\mathcal O}$ omits logarithmic factors in the model parameters and double-logarithmic ones in $(\eps\delta)^{-1}$.

\begin{remark}[Comparison to Previous Work]
Three aspects of our guarantee are noteworthy: the condition number scaling, the dimension scaling, and the dependence on initialization.
We now compare to the previous works \cite{Nawaf26persistence,chen2023}, contributing the basis for the first two aspects, and \cite{BouRabeeOberdoerster2024}, developing the coupling framework whose abstraction in \cite{DPMC} underlies our analysis.

Regarding the condition number aspect, \cite{Nawaf26persistence} established the first complexity guarantee featuring a square-root condition number scaling for Hamiltonian Monte Carlo under general conditions.
In our setting, assuming weaker $L_{H,s}$-spectral Hessian Lipschitzness, the gradient complexity under critical tuning is shown to be
\[ \widetilde{\mathcal O}\Bigl(\bigl((\kappa\Delta)^{1/2}+L_{H,s}^{1/2}m^{-3/4}\Delta^{3/4}\bigl)\log\eps^{-1}\Bigr)\quad\text{with}\quad \Delta\ =\ d+\log(\delta\eps)^{-1} \]
from any $L^2$-warm start such that $\beta_2(\hat\nu)\leq\delta^{-1}$.
Our bounds improve this result in the setting of Frobenius Hessian Lipschitzness and suitable start from $(\kappa d)^{1/2}+L_H^{1/2}m^{-3/4}d^{3/4}$ to $(\kappa^{1/2}+L_H^{1/2}m^{-3/4})d^{1/4}$, up to logarithms.
In the setting of spectral Hessian Lipschitzness and feasible start \cite{Chen19} with $\log\beta_2(\hat\nu)\lesssim d\log\kappa$ however, using $L_H\leq d^{1/2}L_{H,s}$, our result worsens to $(\kappa^{1/2}+L_{H,s}^{1/2}m^{-3/4}d^{1/4})d^{3/4}$ while \cite{Nawaf26persistence} retains the aforementioned guarantee.

Regarding the dimension aspect, \cite{chen2023} first obtained complexity guarantees featuring the $d^{1/4}$ dimension scaling, as predicted by asymptotic optimal-scaling analysis for product targets \cite{Beskos13}, for Hamiltonian Monte Carlo under general conditions.
Specifically, only assuming Cheeger's isoperimetric inequality with coefficient $m^{1/2}$ instead of strong convexity, the gradient complexity of Hamiltonian Monte Carlo in the overdamped regime is shown to be
\[ \widetilde{\mathcal O}\left(\bigl(\kappa+L_H^{2/3}m^{-1}\bigr)\bigl(d+\log(\delta\eps)^{-1}\bigr)^{1/4}\log^{7/4}(\delta\eps)^{-1}\right) \]
from any $L^\infty$-warm start such that $\beta_\infty(\hat\nu)\leq\delta^{-1}$.
Under the additional assumption of strong convexity, our result extends the analysis to critical damping, improving the first factor from $\kappa+L_H^{2/3}m^{-1}$ to $\kappa^{1/2}+L_H^{1/2}m^{-3/4}$.
However, the results of \cite{chen2023} remain valid beyond the strongly convex setting, to which our study is limited.

Regarding the initialization aspect, our study improves upon both \cite{chen2023,Nawaf26persistence}.
On one hand, we extend to general $L^p$-warm starts for $p\in(1,\infty]$ while, on the other hand, introducing new technology that enables deterministic starts in suitable high-probability domains.
The previous work \cite{BouRabeeOberdoerster2024} also considered deterministic starts.
In the setting of \cite{Nawaf26persistence}, the gradient complexity of Metropolis-adjusted kinetic Langevin in the overdamped regime was shown to be
\[ \widetilde{\mathcal O}\left(\kappa^{3/2}\max\bigl((\kappa d)^{1/2},L_{H,s}^{1/2}m^{-3/4}d^{3/4}\bigr)\log^2\eps^{-1}\right) \]
uniformly from start in $\{H\lesssim\kappa d\}$.
The similar dimension scaling of this guarantee to \cite{Nawaf26persistence} originates in the fact that both works localize their analyses to a high-probability sublevel sets of the energy.
In contrast, as \cite{chen2023}, our study localizes to an implicitly given high-probability domain, in which the method is better behaved than in entire energy sublevel sets.
While \cite{chen2023} achieves this localization via conductance machinery, we employ the new coupling technology presented next, refining arguments as used in \cite{BouRabeeOberdoerster2024,Nawaf26persistence}.
\end{remark}

\subsection{High-probability Uniform Mixing via Dobrushin Contraction}\label{sec:mixlocDobHopf}

\Cref{thm:main} for Metropolis-adjusted kinetic Langevin results from a universal framework establishing high-probability uniform mixing of general Markov chains via Dobrushin contraction, which we present now.
Throughout the study of mixing, Dobrushin contraction, that is pointwise total variation contraction, of Markov kernels has been a central concept.
Since Dobrushin's work \cite{Dobrushin56}, global contraction is known to be sufficient for uniform mixing across state space.
In many instances however, global contraction is too strong of an assumption but remains realistic locally.
Harris' ergodic theorems \cite{Harris56} relax contraction from global to local via combination with global Lyapunov drift.
However, due to the Lyapunov condition, these analyses remain inherently global and, in particular, are not capable of capturing fast local equilibration in the presence of slower global mixing \cite{Roberts96b,Qin21}.
A broad class of examples exhibiting this local vs. global mixing at different time scales are Metropolis-adjusted Langevin and Hamiltonian Monte Carlo methods \cite{Diss}.
See the introduction of \cite{DPMC} for a detailed discussion and additional references on the rich history of Harris' theorems.

Several fully local approaches aimed at fast local mixing have been devised, including local variants of conductance \cite{Lovasz93}, spectral gaps \cite{Atchade21}, and Harris' theorem \cite{Yang23}.
In this work, we extend the results in \cite{DPMC}, which use local Dobrushin contraction inside a domain $D$ of state space to conclude geometric mixing up to exit from $D$, a strategy previously described in \cite{Eberle14} for Wasserstein distances.

While different techniques to establish local Dobrushin contraction in suitable $D$ have been and remain to be developed \cite{BouRabeeOberdoerster2024,BouRabeeOberdoersterNUTS1,Diss,Nawaf26persistence}, the stability necessary to control exit from $D$ has remained elusive.
Previous attempts via Lyapunov drift conditions \cite{BouRabeeOberdoerster2024,Yang23}, per-transition exit probabilities \cite{Diss,DPMC}, or restriction to warm initial distributions \cite{Nawaf26persistence} have either required explicit knowledge of $D$ or unnecessarily restricted initialization.
Indeed, on the one hand, domains $D$ in which the desired fast local equilibration holds are usually only implicitly known high-probability regions of the Markov chain's invariant distribution.
On the other hand, there is no reason for stability not to hold from deterministic start in suitable states $x\in D$, i.e., initialization in non-warm distributions $\delta_x$.
The difficulty of showing stability in natural domains has severely limited the ability of localized analyses to reach their full potential.

\Cref{prop:HopfStab} below resolves this theoretical bottleneck by establishing stability in generic high-probability domains $D$ uniformly from start in suitable high-probability subdomains.
It results from a weak-$(1,1)$ maximal inequality that we derive from Hopf's maximal ergodic theorem \cite{Hopf54}.

The following theorem amends the mixing framework based on local Dobrushin contraction from \cite{DPMC} by the stability of \Cref{prop:HopfStab}.
Let $\S$ be a Polish state space endowed with its $\sigma$-algebra $\mathfrak B$, and denote the set of all probability measures on $(\S,\mathfrak B)$ by $\mathcal P(\S)$.

\begin{theorem}[High-probability Uniform Mixing via Dobrushin Contraction]\label{thm:mixlocDobHopf}
Let $\pi$ be a Markov transition kernel on $\S$ with invariant measure $\mu\in\calP(\S)$, and let $\eps,\delta\in(0,1)$.
Suppose there exist $\rho>0$ and a measurable $D\subseteq\S$ such that
\[
    \sup_{x,\tilde x\in D}\TV\bigl(\pi(x,\cdot),\pi(\tilde x,\cdot)\bigr)\ \leq\ 1-\rho\quad\text{and}\quad\mu(D^\comp)\ \leq\ \frac{\eps\,\delta}{4\bigl\lceil\rho^{-1}\log(2/\eps)\bigr\rceil+1}.
\]
Then, there exists a measurable $D'\subseteq D$ such that
\[
    \tmix(\eps,D')\ \leq\ \bigl\lceil\rho^{-1}\log(2/\eps)\bigr\rceil\quad\text{and}\quad\mu(D'^\comp)\ \leq\ \delta.
\]
\end{theorem}

\Cref{thm:mixlocDobHopf} guarantees uniform mixing in a subdomain $D'$ of $\mu$-mass arbitrarily close to unity, assuming Dobrushin contraction to hold in a suitable high-probability domain $D$.
Importantly, only the equilibrium mass of $D$ is needed, which allows its implicit definition.
This considerably simplifies the application of the mixing framework compared to earlier variants.
Notably, the novel stability in high-probability domains underlying this improvement immediately applies to other localized analyses such as \cite{BouRabeeOberdoerster2024,BouRabeeOberdoersterNUTS1,Yang23,Eberle14,Diss,Nawaf26persistence}.

The proof combines \Cref{thm:mixlocDob} and \Cref{prop:HopfStab} and is given in \Cref{sec:pfgeneral}.
A similar result from warm starts is given as \Cref{thm:mixlocDobwarm} in \Cref{sec:warm} below.

\section{Uniform Mixing in High-probability Domains}\label{sec:locDobmix}

We begin this section by recalling the mixing framework via local Dobrushin contraction from \cite{DPMC}.
Subsequently, in \Cref{sec:stabHopf}, we state the novel stability result in high-probability domains and prove it based on Hopf's maximal ergodic theorem.
Combining the two previous results, we then prove the high-probability uniform mixing framework for general Markov chains \Cref{thm:mixlocDobHopf}.
To finish the section, we extend the mixing framework from starts in high-probability domains to $L^p$-warm initial distributions.

\subsection{Local Geometric Mixing via Dobrushin Contraction}\label{sec:mixDob}

Theorem 1 of \cite{DPMC} specializes to the following result, for whose simple proof we refer the reader to \cite{DPMC}.

\begin{theorem}[Local Geometric Mixing via Dobrushin Contraction]\label{thm:mixlocDob}
Let $\pi$ be a Markov transition kernel on $\S$ with invariant measure $\mu\in\calP(\S)$.
Suppose there exist a measurable $D\subseteq\S$ and $\rho>0$ such that
\[
    \TV\bigl(\pi(x,\cdot),\pi(\tilde x,\cdot)\bigr)\ \leq\ 1-\rho\quad\text{for all $x,\tilde x\in D$.}
\]
Then, for all $n\geq0$ and $\nu\in\calP(\S)$, it holds
\[
    \TV(\nu\pi^n,\mu)\ \leq\ \inf_{0\leq k\leq n}\left((1-\rho)^k + \P_{\nu}(T<k) + \P_{\mu}(T<k)\right),
\]
where $T:=\inf\{n\geq0:X_n\notin D\}$ denotes the first exit time from $D$ of the Markov chain with transition kernel $\pi$.
\end{theorem}

Assuming local Dobrushin contraction in a domain $D$, \Cref{thm:mixlocDob} establishes mixing up to exit from $D$.
As discussed in detail in \cite{DPMC}, this enables the theorem to capture fast local equilibration in the presence of slower global mixing.

\subsection{Uniform Stability in High-probability Domains}\label{sec:stabHopf}

Concluding local mixing from \Cref{thm:mixlocDob} requires stability in $D$ both from start in equilibrium $\mu$ as well as in the initial distribution $\nu$.
Domains in which fast local mixing holds commonly arise as high-probability domains with respect to equilibrium.

In such high-probability domains, stability from start in equilibrium holds trivially.
More precisely, if $\pi$ leaves $\mu$ invariant, it holds
\begin{equation}\label{eq:mustab}
    \mathbb P_\mu(T<n)\ \leq\ n\,\mu(D^\comp)\quad\text{for all measurable $D\subseteq\S$ and $n\geq0$.}
\end{equation}
This directly follows from controlling the exit probability in terms of the flow out of $D$ in each transition
\begin{equation}\label{statflow}
    \mathbb P_\mu(T<n)\ \leq\ \mu(D^\comp)\ +\ \sum_{k=0}^{n-2}Q_{\mu\pi^k}(D),
\end{equation}
where $Q_\eta(D)=\int_D\eta(dx)\pi(x,D^\comp)$ for any $\eta\in\mathcal P(\S)$.
By invariance it further holds, for all $k\geq0$,
\[ Q_{\mu\pi^k}(D)\ =\ Q_\mu(D)\ \leq\ \mu\pi(D^\comp)\ =\ \mu(D^\comp). \]

From non-equilibrium starts, stability can be much more involved.
Under restriction to warm starts, stability in suitable high-probability domains is straightforward.
Specifically, by the Markov property and H\"older's inequality, it holds, for all $p\in(1,\infty]$ and $\nu\in\mathcal P(\S)$ such that $\nu\ll\mu$,
\[ \P_\nu(T<n)\ =\ \int_{\S}\P_x(T<n)\frac{d\nu}{d\mu}(x)\,\mu(dx)\ \leq\ \left\|\frac{d\nu}{d\mu}\right\|_{L^p(\mu)}\left(\int_\S\P_x(T<n)^{\frac{p}{p-1}}\,\mu(dx)\right)^{\frac{p-1}{p}}. \]
Hence, by \eqref{eq:mustab},
\begin{equation}\label{eq:stabwarm}
    \P_\nu(T<n)\ \leq\ \beta_p(\nu)\,\P_\mu(T<n)^{\frac{p-1}{p}}\ \leq\ \beta_p(\nu)\bigl(n\mu(D^\comp)\bigr)^{\frac{p-1}{p}}\quad\text{for all $p>1$.}
\end{equation}
Since warm starts arise naturally through the use of analytic techniques on the respective $L^p$ spaces, their direct connection to stability in high-probability domains conveniently yields local mixing via Dobrushin contraction from warm starts, a route we take in \Cref{sec:warm} below.

However, probabilistic coupling approaches are known to accommodate deterministic initialization, for which warmness fails.
Correspondingly, a stability result in high-probability domains that extends beyond warm starts is highly desirable to complement the local mixing framework in \Cref{thm:mixlocDob}.
The following proposition provides such a result.

\begin{proposition}[Uniform Stability in High-probability Domains]\label{prop:HopfStab}
Let $\pi$ be a Markov transition kernel on $\S$ with invariant measure $\mu\in\mathcal P(\S)$, and let $D\subseteq\S$ be a measurable domain.
Then, for all $r>0$, there exists a measurable $D'\subseteq D$ such that
\begin{equation}\label{eq:linex}
    \mathbb P_x(T\leq n)\ \leq\ n\,r\quad\text{for all $x\in D'$ and $n\geq0$,}
\end{equation}
as well as
\[ \mu(D\setminus D')\ \leq\ \frac1r\int_D\mu(dx)\pi(x,D^\comp)\ \leq\ \frac{\mu(D^\comp)}{r}. \]
\end{proposition}

The proposition follows from a weak-$L^1(\mu)$ bound on the maximal function
\[ M(x)\ =\ \ind_D(x)\,\sup_{n\geq1}n^{-1}\P_x(T\leq n), \]
which we show via Hopf's maximal ergodic theorem.
It guarantees stability in any high-probability domain $D$ uniformly from start in an implicitly given subdomain $D'$.
Therein, the equilibrium mass of $D'$ can be fixed arbitrarily close that of $D$.
Notably, stability holds from any initial distribution supported in $D'$, including deterministic initialization in point masses.

\Cref{prop:HopfStab} refines the above average stability estimate under warm start as, by H\"older's inequality,
\begin{equation}\label{eq:hoelderprob}
    \nu(D'^\comp)\ \leq\ \beta_p(\nu)\,\mu(D'^\comp)^{\frac{p-1}{p}}\quad\text{for all measurable $D'\subseteq\S$ and $p\in(1,\infty]$.}
\end{equation}
Hence, a warm start in $X\sim\nu$ satisfies $X\in D'$ with high probability for any $D'$ of sufficiently large equilibrium mass, as provided by the proposition.

\begin{remark}[Comparison with Markov's Inequality]
Under the assumptions of \Cref{prop:HopfStab}, fix $r>0$ and an integer $N\geq1$.
Markov's inequality and \eqref{statflow} imply that
\[
 D_N\ =\ \bigl\{x\in D:\mathbb P_x(T\leq N)\leq Nr\bigr\}
\]
satisfies
\[
 \mu(D\setminus D_N)
 \ \leq\ \frac1{Nr}\int_D\mathbb P_x(T\leq N)\,\mu(dx)
 \ \leq\ \frac1r\int_D\mu(dx)\pi(x,D^\comp).
\]
This fixed-horizon argument already suffices to prove
\Cref{thm:mixlocDobHopf}.
The weak-$(1,1)$ inequality strengthens it by providing a single set $D'$ on which \eqref{eq:linex} holds simultaneously for every horizon, with the same bound on the exceptional mass.
\end{remark}

Turning to the proof, we first state Hopf's maximal ergodic theorem \cite{Hopf54}, see \cite{Garsia65} for a simple proof.

\begin{theorem}[Hopf's Maximal Ergodic Theorem]\label{thm:hopf}
Let $(\S,\mathfrak B,\mu)$ be a measure space and $S$ a non-negative contraction on $L^1(\mu)$, i.e., a linear operator such that $Sf\geq0$ for all non-negative $f\in L^1(\mu)$ and $\|S\|_{L^1(\mu)\to L^1(\mu)}\leq1$.
For $f\in L^1(\mu)$, define $S_0(f):=0$ and $S_n(f):=\sum_{k=0}^{n-1}S^kf$ for $n\geq1$.
Then,
\begin{equation}\label{eq:hopf}
    \int_{\bigl\{\max_{n\leq N}S_n(f)>0\bigr\}}f\,d\mu\ \geq\ 0\quad\text{for all $f\in L^1(\mu)$.}
\end{equation}
If, in addition, $\mu(\S)<\infty$ and $S1\leq1$ almost everywhere, then
\begin{equation}\label{eq:hopfmax}
    \mu\left(\sup_{n\geq1}\frac{S_n(f)}{n}>r\right)\ \leq\ \frac{1}{r}\int_{\S}f\,d\mu\quad\text{for all $r>0$ and non-negative $f\in L^1(\mu)$.}
\end{equation}
\end{theorem}

\begin{proof}[Proof of \Cref{thm:hopf}]
For the proof of \eqref{eq:hopf}, see \cite{Garsia65}.
For the reader's convenience, we recall the classical derivation of \eqref{eq:hopfmax} from \eqref{eq:hopf}.
If $\mu(\S)<\infty$ and $S1\leq1$ almost everywhere, $S^k1\leq1$ for all $k\geq0$ by induction.
Let $r>0$ and $f\in L^1(\mu)$ be non-negative.
As
\[ S_n(f-r)\ =\ S_n(f)-r\sum_{k=0}^{n-1}S^k1\ \geq\ S_n(f)-nr, \]
we have
\[ \left\{\max_{1\leq n\leq N}\frac{S_n(f)}{n}>r\right\}\ \subseteq\ \left\{\max_{n\leq N}S_n(f-r)>0\right\}\quad\text{for all $N\geq1$.} \]
By \eqref{eq:hopf} and $f\geq0$,
\[ 0\ \leq\ \int_{\bigl\{\max_{n\leq N}S_n(f-r)>0\bigr\}}(f-r)\,d\mu\ \leq\ \int_{\S}f\,d\mu-r\,\mu\left(\max_{n\leq N}S_n(f-r)>0\right), \]
which yields
\[ \mu\left(\sup_{n\geq1}\frac{S_n(f)}{n}>r\right)\ \leq\ \lim_{N\to\infty}\mu\left(\max_{n\leq N}S_n(f-r)>0\right)\ \leq\ \frac{1}{r}\int_{\S}f\,d\mu. \]
\end{proof}

\begin{proof}[Proof of \Cref{prop:HopfStab}]
The linear operator
\[ Sf(x)\ :=\ \ind_D(x)\int_D\pi(x,dy)f(y) \]
defines a non-negative contraction on $L^1(\mu)$.
Indeed, by invariance,
\[ \int_\S|Sf|\,d\mu\ \leq\ \int_\S|f|\,d(\mu\pi)\ =\ \int_\S|f|\,d\mu\quad\text{for all $f\in L^1(\mu)$.} \]
Moreover, $S1(x)=1_D(x)\pi(x,D)\leq1$.
Let $f(x):=\ind_D(x)\,\pi(x,D^\comp)\in L^1(\mu)$ and define
\[ g_n(x)\ :=\ \ind_D(x)\,\P_x(T\leq n)\quad\text{for $n\geq0$.} \]
By conditioning on $X_1$ and using the Markov property, we see the recursion
\begin{align*}
    g_n(x)\ &=\ \P_x\bigl(T\leq n,X_1\in D\bigr)+\P_x(X_1\in D^\comp)\ =\ \int_D g_{n-1}(y)\,\pi(x,dy)+f(x)\\
    &=\ Sg_{n-1}(x)+f(x)\quad\text{for $n\geq1$ and $x\in D$.}
\end{align*}
Therefore $g_n=Sg_{n-1}+f$ for all $n\geq1$, iterating which shows
\[ g_n\ =\ \sum_{k=0}^{n-1}S^kf\ =\ S_n(f). \]
Let $r>0$ and define
\[ D'\ :=\ \bigl\{x\in D:g_n(x)\leq nr\text{ for all $n\geq0$}\bigr\}. \]
Then, by definition,
\[ \mathbb P_x(T\leq n)\ \leq\ nr\quad\text{for all $x\in D'$ and $n\geq0$,} \]
as well as, by \eqref{eq:hopfmax},
\[ \mu(D\setminus D')\ \leq\ \mu\left(\sup_{n\geq1}\frac{g_n}{n}>r\right)\ \leq\ \frac1r\int_D\mu(dx)\pi(x,D^\comp)\ \leq\ \frac{\mu(D^\comp)}r, \]
where the last step uses invariance of $\mu$ under $\pi$.
\end{proof}

\subsection{High-probability Uniform Mixing via Dobrushin Contraction}\label{sec:pfgeneral}

We are now positioned to prove \Cref{thm:mixlocDobHopf} by combination of local geometric mixing via Dobrushin contraction from \Cref{thm:mixlocDob} with the uniform stability in high-probability domains of \Cref{prop:HopfStab}.

\begin{proof}[Proof of \Cref{thm:mixlocDobHopf}]
Let $\pi$ be a Markov transition kernel on $\S$ with invariant measure $\mu\in\calP(\S)$, and let $\eps,\delta\in(0,1)$.
Further, let $\rho>0$, $n=\bigl\lceil\rho^{-1}\log(2/\eps)\bigr\rceil$, and $D\subseteq\S$ be measurable such that
\[ \sup_{x,\tilde x\in D}\TV\bigl(\pi(x,\cdot),\pi(\tilde x,\cdot)\bigr)\ \leq\ 1-\rho\quad\text{and}\quad\mu(D^\comp)\ \leq\ \frac{\eps\,\delta}{2n(1+\delta)+\eps}. \]
The second part follows from $\mu(D^\comp)\leq\frac{\eps\delta}{4n+1}$ as assumed in the theorem.
We can assume $\mu(D^\comp)>0$ as otherwise, after removing a $\mu$-null set from $D$, the chain remains in $D$ almost surely, in which case the local Dobrushin contraction implies geometric mixing with rate $\rho$ uniformly from start in $D$.
In addition, the second part of the previous display further implies $\mu(D^\comp)<\delta$ so that we may fix $r=\frac{\mu(D^\comp)}{\delta-\mu(D^\comp)}>0$.
By \Cref{prop:HopfStab}, there exists a measurable $D'\subseteq D$ such that, for all $n\geq0$
\[ \sup_{x\in D'}\mathbb P_x(T\leq n)\ \leq\ n\,r,\quad\text{and}\quad\mu(D\setminus D')\ \leq\ \frac{\mu(D^\comp)}{r}. \]
The second inequality is equivalent to $\mu(D'^\comp)\leq(1+r^{-1})\mu(D^\comp)=\delta$ by our choice of $r$.
By \Cref{thm:mixlocDob}, the first inequality in the previous display and \eqref{eq:mustab},
\[ \sup_{x\in D'}\TV\bigl(\pi^n(x,\cdot),\mu\bigr)\ \leq\ (1-\rho)^n\ +\ n\,\bigl(\mu(D^\comp)+r\bigr)\quad\text{for all $n\geq0$.} \]
Fixing $n=\bigl\lceil\rho^{-1}\log(2/\eps)\bigr\rceil$ bounds the first term above by $\eps/2$.
The second term satisfies
\[ n\,\bigl(\mu(D^\comp)+r\bigr)\ =\ \frac{n\mu(D^\comp)(1+\delta-\mu(D^\comp))}{\delta-\mu(D^\comp)}\ \leq\ \frac{n\mu(D^\comp)(1+\delta)}{\delta-\mu(D^\comp)}\ \leq\ \frac\eps2,  \]
where we inserted $r$ and used that our assumption on $\mu(D^\comp)$ reads $\mu(D^\comp)\leq\frac{\eps\delta}{2n(1+\delta)+\eps}$, which is equivalent to the last inequality in the previous display.
Thus $\TV\bigl(\pi^n(x,\cdot),\mu\bigr)\leq\eps$ for all $x\in D'$, which implies $\tmix(\eps,D')\leq n$, finishing the proof.
\end{proof}

\subsection{Geometric Mixing via Dobrushin Contraction from Warm Starts}\label{sec:warm}

To finish the section, we make the mixing framework via local Dobrushin contraction accessible from warm starts.
Therefore, we combine \Cref{thm:mixlocDob} with stability from warm initial distributions shown in \eqref{eq:stabwarm}.

\begin{theorem}[Geometric Mixing via Dobrushin Contraction from Warm Starts]\label{thm:mixlocDobwarm}
Let $\pi$ be a Markov transition kernel on $\S$ with invariant measure $\mu\in\calP(\S)$, and let $\eps\in(0,1)$, $p\in(1,\infty]$, and $\nu\in\mathcal P(\S)$ be $L^p$-warm.
Suppose there exist $\rho>0$ and a measurable $D\subseteq\S$ such that
\begin{equation}\label{eq:asswarm}
    \sup_{x,\tilde x\in D}\TV\bigl(\pi(x,\cdot),\pi(\tilde x,\cdot)\bigr)\ \leq\ 1-\rho\quad\text{and}\quad\mu(D^\comp)\ \leq\ \frac{1}{\bigl\lceil\rho^{-1}\log(2/\eps)\bigr\rceil}\left(\frac{\eps}{4\beta_p(\nu)}\right)^{\frac{p}{p-1}}.
\end{equation}
Then, 
\[ \tmix(\eps,\nu)\ \leq\ \bigl\lceil\rho^{-1}\log(2/\eps)\bigr\rceil. \]
\end{theorem}

The theorem guarantees mixing from warm start, given that Dobrushin contraction holds in a domain $D$ of sufficiently large equilibrium mass.
Alternatively to combining \Cref{thm:mixlocDob} directly with stability from warm starts, one may obtain a similar result as a corollary of the high-probability uniform mixing framework \Cref{thm:mixlocDobHopf}.
Therefore, let $D'\subseteq\S$ be a suitable high-probability subdomain of $D$.
By the Markov property,
\[ \TV\bigl(\nu\pi^n,\mu\bigr)\ \leq\ \sup_{x\in D'}\TV\bigl(\pi^n(x,\cdot),\mu\bigr)\ +\ \nu(D'^\comp) \]
and, by \eqref{eq:hoelderprob}, a warm start $X\sim\nu$ is likely contained in $D'$.
Thus, we may control the first term via \Cref{thm:mixlocDobHopf} while controlling the second via \eqref{eq:hoelderprob}.
This detour however requires $\mu(D^\comp)$ to be smaller by a factor $\eps$ compared to \eqref{eq:asswarm}.

\begin{proof}[Proof of \Cref{thm:mixlocDobwarm}]
Assume the setting stated in the theorem and let $n=\bigl\lceil\rho^{-1}\log(2/\eps)\bigr\rceil$.
By \Cref{thm:mixlocDob} together with stability from equilibrium \eqref{eq:mustab} and warm start \eqref{eq:stabwarm}, it holds
\[ \TV\bigl(\nu\pi^n,\mu\bigr)\ \leq\ (1-\rho)^n\ +\ n\mu(D^\comp)\ +\ \beta_p(\nu)\bigl(n\mu(D^\comp)\bigr)^{\frac{p-1}{p}}. \]
Our choice of $n$ implies the first term on the right hand side to be $\leq\eps/2$.
The other two terms are each $\leq\eps/4$ by the assumption on $\mu(D^\comp)$.
Specifically, since $\beta_p(\nu)\geq1$ and $\frac{p}{p-1}\geq1$,
\[ \mu(D^\comp)\ \leq\ \frac{1}{n}\left(\frac{\eps}{4\beta_p(\nu)}\right)^{\frac{p}{p-1}}\ \leq\ \frac{\eps}{4n}. \]
Hence, $\TV\bigl(\nu\pi^n,\mu\bigr)\leq\eps$, so that $\tmix(\eps,\nu)\leq n$ as claimed.
\end{proof}

\section{Mixing of Metropolis-adjusted Kinetic Langevin}\label{sec:proofs}

Our result \Cref{thm:main} for Metropolis-adjusted kinetic Langevin with persistent uniforms is an application of the general framework \Cref{thm:mixlocDobHopf}.
Therefore, we require \emph{Dobrushin contraction of Metropolis-adjusted kinetic Langevin in a high-probability domain}, which we establish in three steps:
\begin{enumerate}[(i)]
\item
\Cref{prop:obabo-dobrushin} shows Dobrushin contraction in a high-probability domain over a suitably long epoch of $E$ unadjusted OBABO steps,
\item
\Cref{prop:pu-acceptance} shows the probability of a rejection occurring during such an epoch of the Metropolis-adjusted method to be controlled uniformly in a high-probability domain.
\item
\Cref{lem:MAKLApuDob} combines (i) and (ii) into Dobrushin contraction of the Metropolis-adjusted method in a high-probability domain.
It is based on the \emph{high-acceptance} idea, which reduces Dobrushin contraction of a Metropolis-adjusted Markov chain over an epoch to the contraction of its underlying unadjusted chain, given that the probability of rejection over the epoch is controlled.
Specifically, it holds
\[ \TV\bigl(\Pi(z,\cdot),\Pi(\tilde z,\cdot)\bigr)\ \leq\ \TV\bigl(\Pi_u(z,\cdot),\Pi_u(\tilde z,\cdot)\bigr)\ +\ 2\,\sup_{z\in D}\mathbb P_z\left(\left(\bigcap\nolimits_{1\leq k\leq E}A_k\right)^\comp\right) \]
for all $z,\tilde z\in D$ in any measurable domain $D\subseteq\R^{2d}$.
\end{enumerate}
\noindent
The high-acceptance idea was first proposed in \cite{BouRabeeOberdoerster2024}, see \cite{Diss} for a systematic development, and has since been applied to kinetic Langevin \cite{BouRabeeOberdoerster2024,Nawaf26persistence}, Hamiltonian Monte Carlo \cite{Diss}, and the No-U-turn sampler \cite{BouRabeeOberdoersterNUTS1,NUTS2}.

We now state the propositions making up steps (i) and (ii) as well as the lemma combining them in (iii).
Their proofs and the proof of \Cref{thm:main} are given below.

\begin{proposition}[High-Probability Unadjusted Dobrushin Contraction over Epoch]\label{prop:obabo-dobrushin}
Suppose $\mu$ satisfies \Cref{ass:1}, and, as in \Cref{def:MAKLApu}, let $\Pi_u$ be the unadjusted transition kernel of an epoch of
$E\in\mathbb N$ steps of OBABO with friction $\gamma=2\sqrt m$ and step size $h>0$.
For $R\geq0$, set
\[
 D_R\ =\ \bigl\{z\in\R^{2d}:H(z)\leq R\bigr\}\quad\text{and}\quad\Lambda_R\ =\ \log\max\bigl(e,\kappa,L_H,d,R\bigr).
\]
There exists
\begin{equation}\label{eq:obabo-parameters}
 \bar h^{-1}\ \asymp\ \bigl(L^{1/2}+L_H^{1/2}m^{-1/4}\bigr)(d+R)^{1/4}\Lambda_R^{1/4}
\end{equation}
such that, for all $h\leq\bar h$, there exists $E\asymp\Lambda_R/(h\sqrt m)$, for which
\begin{equation}\label{eq:obabo-dobrushin}
 \TV\bigl(\Pi_u(z,\cdot),\Pi_u(\tilde z,\cdot)\bigr)\ \leq\ \frac12\quad\text{for all $z,\tilde z\in D_R$.}
\end{equation}
\end{proposition}

Note that, as shown in \cite[Lem.~5.4]{Nawaf26persistence}, $D_R$ for suitable $R$ is a high-probability domain with respect to equilibrium.

\begin{proposition}[High-probability High-acceptance over Epoch]\label{prop:pu-acceptance}
Suppose $\mu$ satisfies \Cref{ass:1} and, as in \Cref{def:MAKLApu}, let $(A_k)_{1\leq k\leq E}$ denote the acceptance events during an epoch of $E\in\mathbb N$ steps of Metropolis-adjusted kinetic Langevin with persistent uniforms with friction $\gamma\geq0$ and step size $h>0$.
Let $C,R\geq0$ and $\delta\in(0,1/2)$.
With $\Lambda_R$ as in \Cref{prop:obabo-dobrushin}, there exists
\begin{equation}\label{eq:pu-step-size}
 \bar h^{-1}\ \asymp_C\ \bigl(L^{1/2}+\Lambda_R^{1/2}L_H^{1/2}m^{-1/4}\bigr)\bigl(d+\log(E/\delta)\bigr)^{1/4}\log^{1/2}(E/\delta),
\end{equation}
such that, for all $h\leq\bar h$ satisfying $E\leq C\Lambda_R/(h\sqrt m)$, there is a measurable $D\subseteq\mathbb R^{2d}$ with
\begin{equation}\label{eq:pu-good-domain}
 \hat\mu(D^\comp)\ \leq\ \delta\quad\text{and}\quad\sup_{z\in D}\mathbb P_z\left(
   \left(\bigcap\nolimits_{1\leq k\leq E}A_k\right)^\comp\right)\ \leq\ \frac18.
\end{equation}
\end{proposition}

Combining \Cref{prop:obabo-dobrushin,prop:pu-acceptance} via the high-acceptance idea yields high-probability adjusted Dobrushin contraction.

\begin{lemma}[High-probability Adjusted Dobrushin Contraction over Epoch]\label{lem:MAKLApuDob}
As in \Cref{def:MAKLApu}, let $\Pi$ be the transition kernel of an epoch of $E\in\mathbb N$ steps of Metropolis-adjusted kinetic Langevin with persistent uniforms.
For any measurable domain $D\subseteq\R^{2d}$, it holds, for all $z,\tilde z\in D$,
\begin{equation}\label{eq:puDobsep}
    \TV\bigl(\Pi(z,\cdot),\Pi(\tilde z,\cdot)\bigr)\ \leq\ \TV\bigl(\Pi_u(z,\cdot),\Pi_u(\tilde z,\cdot)\bigr)\ +\ 2\,\sup_{z\in D}\mathbb P_z\left(\left(\bigcap\nolimits_{1\leq k\leq E}A_k\right)^\comp\right).
\end{equation}
Further, suppose $\mu$ satisfies \Cref{ass:1} and let $\gamma=2\sqrt m$, $\delta\in(0,1/2)$.
There exists a measurable domain $D\subseteq\R^{2d}$, a step size $h>0$, and an epoch length
\[ E\ \lesssim\ \bigl(\kappa^{1/2}\Lambda+L_H^{1/2}m^{-3/4}\Lambda^{3/2}\bigr)\bigl(d+\log\delta^{-1}+\Lambda\bigr)^{1/4}\bigl(\log\delta^{-1}+\Lambda\bigr)^{1/2}, \]
where $\Lambda=\log\max\bigl(e,\kappa,L_H,L_H/m^{3/2},d,\log\delta^{-1}\bigr)$, such that
\[ \sup_{z,\tilde z\in D}\TV\bigl(\Pi(z,\cdot),\Pi(\tilde z,\cdot)\bigr)\ \leq\ \frac34\quad\text{and}\quad\hat\mu(D^\comp)\ \leq\ \delta. \]
\end{lemma}

\subsection{High-probability Unadjusted Dobrushin Contraction over Epoch}\label{sec:obabo-contraction-proof}

We start with the proof of \Cref{prop:obabo-dobrushin}.
For each starting point, we approximate the distribution after one epoch in total variation by a distribution with small relative entropy with respect to $\hat\mu$.
Pinsker's inequality then compares two starts.
We follow the entropy argument in \cite[Sec.~4--6]{Nawaf26persistence}, which adds the mean inner product of velocity with an optimal transport displacement \cite{Lu26entropy} and subtracts the gradient-square term from the leapfrog energy identity \cite{Bussi07}.
Here we center the latter term and estimate expectations of the errors before taking absolute values.
The resulting dimensional term in the one-step error is of order $h^5(L^2+L_H^2)d$.
This uses Hessian Lipschitzness in Frobenius norm, a stronger assumption than the spectral norm condition in \cite{Nawaf26persistence}.
Bounds on energy over elapsed time and on entropy after a fixed elapsed time make the conclusion uniform over all $h\leq\bar h$.

We state the principal estimate in the normalization $m=1$, $\gamma=2$.
Write $\varphi=\mathcal N(0,I_d)$.
For $p\in\mathcal P(\mathbb R^{2d})$ with $\mathcal H(p)=\KL(p\mid\hat\mu)<\infty$,
let $p_x$ be its position marginal and $T_p$ the optimal transport map
from $p_x$ to $\mu$ for squared Euclidean cost.
For a sufficiently small universal $\eta>0$, define
\begin{equation}\label{eq:readable-centered-functional}
 \begin{aligned}
 \mathcal C(p)\ &=\ \int v\cdot\bigl(x-T_p(x)\bigr)\,p(dx,dv),
 &G_h(x)\ &=\ \tfrac18h^2\bigl|\nabla U(x)\bigr|^2,\\
 \mathcal F_h(p)\ &=\ \mathcal H(p)+\eta\mathcal C(p)-p(G_h)+\mu(G_h).
 \end{aligned}
\end{equation}
Here $p(f)=\int f\,dp$. Finite relative entropy implies finite second
moments, so all terms are defined. Let $O$ be an OU half-step and
$K_h=O\Phi_hO$, with $pK_h=\bigl((pO)\Phi_h\bigr)O$.
Deterministic maps act by pushforward. For $(X,v_-)\sim pO$, put
$V=v_--(h/2)\nabla U(X)$ and define
\[
 T_r(p)\ =\ \E\bigl[|V|^2\mathbf1_{\{|V|>r\}}\bigr],\qquad
 t_r\ =\ \varphi\bigl(|v|^2\mathbf1_{\{|v|>r\}}\bigr).
\]
Thus $T_r$ includes the first OU half-step and first kick.
Implicit constants and smallness thresholds are universal unless a dependence is indicated.

\begin{lemma}[One-step Entropy Decrease]\label{lem:readable-local}
Suppose $\mu$ satisfies \Cref{ass:1} with $m=1$, and take $\gamma=2$.
There is a universal $c>0$ such that the following holds for every $p\in\mathcal P(\mathbb R^{2d})$ with $\mathcal H(p)<\infty$.
If $h>0$, $r\geq2\sqrt d$, and
\begin{equation}\label{eq:unified-local-cap}
 h^4\bigl((L^2+L_H^2)d+L_H^2r^2\bigr)\quad\text{is sufficiently small,}
\end{equation}
then
\begin{align}
 &\max\bigl(\tfrac12\mathcal H(p)-\mathcal F_h(p),
 \mathcal F_h(p)-2\mathcal H(p)\bigr)\ \lesssim\ h^4L^2d,
 \label{eq:unified-functional-comparison}\\
 &\mathcal F_h(pK_h)-\bigl(1-ch\bigr)\mathcal F_h(p)
\ \lesssim\ h^5(L^2+L_H^2)d
       +\bigl(Lh^2+L^2h^3\bigr)T_r(p)+Lh^3t_r.
 \label{eq:unified-local-recursion}
\end{align}
\end{lemma}

The proof of \Cref{lem:readable-local} is given in \Cref{app:A}. Centering $G_h$
bounds $|p(G_h)-\mu(G_h)|$ by a small multiple of entropy plus
$h^4L^2d$. Estimating the expectations of the leapfrog errors leaves the
$h^5(L^2+L_H^2)d$ term in \eqref{eq:unified-local-recursion}. Its contribution
after iteration is
\begin{equation}\label{eq:structure-accumulated-error}
 h^5(L^2+L_H^2)d\sum\nolimits_{k\geq0}(1-ch)^k\ =\ c^{-1}h^4(L^2+L_H^2)d.
\end{equation}
The corresponding error terms in \cite[Thm.~2]{Nawaf26persistence}
are of order $h^2$ times energy and velocity moments. Our fourth-order
error, together with the velocity cutoff in
\eqref{eq:unified-local-cap}, accounts for the fourth root in the
step-size bound. The remaining task is to initialize the entropy and
control the velocity tails without losing uniformity as $h$ decreases.

\begin{proof}[Proof of \Cref{prop:obabo-dobrushin}]
First normalize $m=1$ and $\gamma=2$ by setting
$x'=\sqrt m\,x$, $U'(x')=U(x'/\sqrt m)$ and $h'=h\sqrt m$.
This replaces $L,L_H$ by $L/m,L_H/m^{3/2}$ and preserves energy,
relative entropy and total variation. Suppress primes, but retain
$\Lambda_R$ from the original model. Its only uses are
\begin{equation}\label{eq:readable-log-controls}
 \Lambda_R\ \geq\ 1,\qquad \log L\ \leq\ \Lambda_R,\qquad \log(R+d)\ \leq\ 2\Lambda_R.
\end{equation}

Two estimates from \Cref{app:A} supply the initial entropy and the tails.
\Cref{lem:readable-smoothing} bounds entropy after elapsed time $t=jh$
with $1\leq t\leq2$ when $Lh\leq1$, replacing the $\log(1/h)$ term in the one-step
bound \cite[Lem.~4.13]{Nawaf26persistence}. When $Lh>1$, that bound
already suffices because $\log(1/h)\leq\log L$.
\Cref{lem:readable-localization} controls large energy and intermediate
velocity by stopping the estimate of \cite[Lem.~5.5]{Nawaf26persistence},
in place of the union over steps in \cite[Lem.~5.6]{Nawaf26persistence}.
For $Lh>1$, we use a comparison chain that takes one unadjusted step
and then transitions with $K_h^\#$. In these transitions, leapfrog is
replaced by velocity reversal if the intermediate velocity exceeds $v_*$.
Driven by the same Gaussian variables, the two chains agree outside the exceptional event. The comparison chain
has controlled energy moments and satisfies a one-sided comparison
of $\mathcal F_h$ with an unadjusted update.

For sufficiently large universal $K,A$ and sufficiently small $c_A>0$, set
\begin{equation}\label{eq:readable-epoch-choice}
 \bar h\ =\ \frac{c_A}{\bigl(\sqrt L+\sqrt{L_H}\bigr)\bigl((R+d)\Lambda_R\bigr)^{1/4}},\qquad
 E\ =\ \left\lceil\frac{A\Lambda_R}{h}\right\rceil,\qquad T\ =\ (A+3)\Lambda_R.
\end{equation}
Take any $0<h\leq\bar h$, choose $\epsilon=e^{-K(T+\Lambda_R)}$ and
$R_*\asymp_A(R+d)\Lambda_R$, and use $r=3\sqrt{R_*}$, $v_*=2r$.
The radius constant and then $c_A$ can be chosen to satisfy all
conditions of the three lemmas, as verified in
\Cref{sec:obabo-epoch-parameters}.
For $z\in D_R$ and all $0\leq k\leq N$, let
\begin{align*}
 &Lh \leq 1:\ p_k\ =\ \delta_zK_h^{j+k}\quad\text{with $j=\lceil1/h\rceil$ and $N=E-j$,}\\
 &Lh > 1:\ p_k\ =\ (\delta_zK_h)(K_h^\#)^k\quad\text{with $N=E-1$.}
\end{align*}
The starting-entropy estimate implies
$\mathcal H(p_0)\lesssim L(R+d)\Lambda_R$, and $Nh\geq Eh-2$.
The velocity tails in \Cref{lem:readable-localization} bound the error
in \eqref{eq:unified-local-recursion} by a universal multiple of
$h\bigl[h^4(L^2+L_H^2)d+P\sqrt\epsilon\bigr]$, where
$\log(1+P)\lesssim T+\Lambda_R$.
For the comparison chain, \eqref{eq:unified-transfer} has the same
bound, since $h^{-1}<L$. Thus iteration of the local estimate, followed
by \eqref{eq:unified-functional-comparison}, shows that $p_N$ satisfies
\[
 \begin{aligned}
 \TV\bigl(\delta_zK_h^E,p_N\bigr)\ &\leq\ \epsilon,\\
 \KL(p_N\mid\hat\mu)\ &\lesssim\
 e^{-c(Eh-2)}L(R+d)\Lambda_R+h^4(L^2+L_H^2)d+P\sqrt\epsilon.
 \end{aligned}
\]
The first distance is zero when $Lh\leq1$. The complete iteration and
an explicit admissible $P$ are recorded in \Cref{sec:obabo-epoch-parameters}.
Choose $K$ large for the tail term, $A$ large for the transient term,
then the radius constant, and finally $c_A$ small for the step-size
conditions and the discretization term. These choices ensure
$\epsilon\leq1/16$ and $\KL(p_N\mid\hat\mu)\leq1/32$.
Pinsker's inequality implies
\[
 \TV\bigl(\delta_zK_h^E,\hat\mu\bigr)\ \leq\ 1/16+\sqrt{(1/32)/2}\ =\ 3/16.
\]
The triangle inequality proves the desired bound for two starts.

Restoring the original parameters, \eqref{eq:readable-epoch-choice} becomes \eqref{eq:obabo-parameters}, and $E\asymp\Lambda_R/(h\sqrt m)$.
\end{proof}

\subsection{High-probability High-acceptance over Epoch}\label{sec:obabo-acceptance-proof}

Next, we prove \Cref{prop:pu-acceptance}.
The persistent-uniform criterion in \cite[Prop.~5.1]{Nawaf26persistence} reduces acceptance of the whole epoch to a bound on the partial sums of leapfrog energy errors.
The proof in \cite[Lem.~3 and~8, Sec.~5.4]{chen2023} estimates the full energy error through moments under $\hat\mu$ and compares the densities of the numerical states to the density of $\hat\mu$, restricted to paths with small preceding errors.
Here we subtract the gradient-square endpoint difference and estimate the remaining error with the weight $e^{sX}$ defined below.
Summing differences of $p$-th roots of densities then accumulates only the remainder proportional to $L_H$.
We include each Gaussian refresh variable in the state so that the OU steps become volume-preserving rotations.
The weighted estimates use strong convexity, which is not required by \cite{chen2023}.

The proof first estimates the partial sums for an initial point distributed according to $\hat\mu$, then conditions on that point to obtain $D$.
Write $\varphi=\mathcal N(0,I_d)$ and set $F=|\nabla U|^2$, $X=F-\mu(F)$, $s=h^2/8$, and
$\mathcal R=\Delta H_h-s\bigl(F\circ\Phi_h-F\bigr)$, regarding $F$ as a function of position on phase space.
Thus $sF=G_h$.

\begin{lemma}[Corrected Leapfrog Energy Error]\label{lem:acceptance-local}
Suppose $\mu$ satisfies \Cref{ass:1}.
If $h>0$ and $h^2L\sqrt d$ is sufficiently small, then,
for every integer $p\geq2$,
\begin{equation*}
 \bigl\|e^{sX/p}-1\bigr\|_{L^p(\hat\mu)}\ \lesssim\ h^2L\sqrt{d+p}\quad\text{and}\quad\left(\int e^{sX}|\mathcal R|^p\,d\hat\mu\right)^{1/p}\ \lesssim\ h^3L_Hp\sqrt{d+p}.
\end{equation*}
\end{lemma}

The first estimate bounds the $L^p(dz)$ distance between the $p$-th roots of the densities of $\hat\mu$ and the unnormalized measure $e^{sX}\hat\mu$, where $dz$ is Lebesgue measure on phase space. The second bounds the corrected leapfrog error under the latter measure.
The proof of \Cref{lem:acceptance-local} is given in \Cref{app:B}.

\begin{proof}[Proof of \Cref{prop:pu-acceptance}]
Suppose $\mu$ satisfies \Cref{ass:1}, let $E\geq1$ be an integer, and let $h>0$ and $\gamma\geq0$.
First assume only that $h^2L\sqrt d$ is sufficiently small. We derive an estimate for every integer $p\geq2$ and every $b>0$, then impose the duration and step-size conditions of the proposition.

For the acceptance criterion, let $Z^u_0=z$ and $Z^u_k=\Phi_{h,\gamma}^{\xi_k,\xi'_k}(Z^u_{k-1})$ be the unadjusted path from \Cref{def:MAKLApu}.
Define the partial sums of its leapfrog energy errors by
\[ W_j\ =\ \sum_{k=1}^j\Delta H_h\!\left(O_{\gamma,h/2}^{\xi_k}(Z^u_{k-1})\right)\quad\text{for all $0\leq j\leq E$, and}\quad M_E\ =\ \max_{0\leq j\leq E}W_j. \]
Here $W_0=0$. On acceptance of the first $j-1$ steps, $U_{j-1}=U_0e^{W_{j-1}}$, so the next step is accepted if $U_0e^{W_j}<1$.
Conditional on the Gaussian noises, all steps are therefore accepted with probability $e^{-M_E}$.
Hence
\begin{equation}\label{eq:acceptance-threshold}
    \mathbb P_z\left(\bigcap\nolimits_{k=1}^E A_k\right)\ =\ \mathbb E_z\,e^{-M_E}.
\end{equation}

Thus it suffices to bound the positive maximum $M_E$, rather than just the terminal sum $W_E$.
We use Lebesgue $L^p$ norms of differences of $p$-th roots of densities, since volume-preserving changes of variables leave these norms unchanged.
Time reversibility of leapfrog converts the weighted local moment into such a difference for one leapfrog step. Including the OU noise coordinates allows us to sum these differences over the unadjusted path without assuming its intermediate distributions are stationary.
For the moment estimates, let $Z_0^u\sim\hat\mu$ independently of the noises. Fix an integer $p\geq2$ and, on phase space with Lebesgue measure $dz$, set
\[ f\ =\ \left(\frac{d\hat\mu}{dz}\right)^{1/p}\quad\text{and}\quad g\ =\ f\,e^{sX/p}. \]
Then $\|f-g\|_{L^p(dz)}$ is the first norm in \Cref{lem:acceptance-local}, and $g\circ\Phi_h=g\,e^{-\mathcal R/p}$.
The leapfrog map preserves volume and satisfies $\Phi_hS\Phi_h=S$ \cite{Nawaf18}.
Thus $I=S\Phi_h$ is a volume-preserving involution with $\mathcal R\circ I=-\mathcal R$ and $g\circ I=g\,e^{-\mathcal R/p}$.
Changing variables by $I$ on $\{\mathcal R<0\}$ yields
\begin{align*}
 \bigl\|g\circ\Phi_h-g\bigr\|_{L^p(dz)}^p
\ &=\ 2\int_{\{\mathcal R>0\}}
       e^{sX}\bigl(1-e^{-\mathcal R/p}\bigr)^p\,d\hat\mu\ \lesssim\ \frac1{p^p}\int e^{sX}|\mathcal R|^p\,d\hat\mu.
\end{align*}
After taking $p$-th roots, \Cref{lem:acceptance-local} bounds the increment by
\begin{equation}\label{eq:readable-root-increment}
    \bigl\|g\circ\Phi_h-g\bigr\|_{L^p(dz)}\ \lesssim\ h^3L_H\sqrt{d+p}.
\end{equation}
To apply the same change-of-variables argument through the OU steps, retain their Gaussian inputs and consider the extended space
\[ \zeta\ =\ \bigl(x,v,\eta_1,\eta'_1,\ldots,\eta_E,\eta'_E\bigr)\ \in\ \mathbb R^{2d}\times\bigl(\mathbb R^d\times\mathbb R^d\bigr)^E\quad\text{with}\quad\bar\mu\ =\ \hat\mu\otimes\varphi^{\otimes2E}. \]
Initially $(\eta_k,\eta'_k)=(\xi_k,\xi'_k)$, independently of $(x,v)\sim\hat\mu$.
Relative to Lebesgue measure $d\zeta$, set
\[
 \bar f\ =\ \left(\frac{d\bar\mu}{d\zeta}\right)^{1/p}\quad\text{and}\quad\bar g\ =\ \bar f\,e^{sX/p}.
\]
Let $\bar\Phi_h$ act as $\Phi_h$ on the phase coordinates and fix the noises. Let $\mathsf O_k$ rotate $(v,\eta_k)$ by
\[
 (v,\eta_k)\ \longmapsto\ \left(e^{-\gamma h/2}v+\sqrt{1-e^{-\gamma h}}\eta_k,-\sqrt{1-e^{-\gamma h}} v+e^{-\gamma h/2}\eta_k\right),
\]
leaving all other coordinates unchanged, and define $\mathsf O'_k$ in the same way with $\eta'_k$.
These rotations preserve volume and both $\bar f,\bar g$. Set
\[
 \Psi_0\ =\ \mathrm{id}\quad\text{and}\quad
 \Psi_j\ =\ \mathsf O'_j\circ\bar\Phi_h\circ\mathsf O_j\circ\Psi_{j-1}.
\]
Each noise coordinate retains its initial value until it is used in a refresh.
Its updated value then remains in the extended state but is not used in any later refresh.
Thus the phase space projection of $\Psi_j$ is exactly the unadjusted path.
Regarding $W_j$ as a function of the initial extended point, each leapfrog step multiplies $\bar f$ by $e^{-\Delta H_h/p}$, with $\Delta H_h$ evaluated at the phase input, and each rotation leaves $\bar f$ unchanged.
Thus, by induction,
\[
 \bar f\circ\Psi_j\ =\ \bigl(\bar f\circ\Psi_{j-1}\bigr)e^{-(W_j-W_{j-1})/p}\ =\ \bar f\,e^{-W_j/p}.
\]
By invariance under the rotations,
\[
 \bar g\circ\Psi_j-\bar g\circ\Psi_{j-1}
\ =\ \bigl(\bar g\circ\bar\Phi_h-\bar g\bigr)
       \circ\mathsf O_j\circ\Psi_{j-1}.
\]
Changing variables by the volume-preserving maps preserves the Lebesgue $L^p$-norm. Insert $\bar g\circ\Psi_j$ and $\bar g$ between $\bar f\circ\Psi_j$ and $\bar f$, then sum the displayed differences. Integrating the Gaussian factors then yields
\begin{align*}
 \bigl\|e^{-W_j/p}-1\bigr\|_{L^p(\bar\mu)}
\ &=\ \bigl\|\bar f\circ\Psi_j-\bar f\bigr\|_{L^p(d\zeta)}
\ \leq\ 2\|f-g\|_{L^p(dz)}
          +j\bigl\|g\circ\Phi_h-g\bigr\|_{L^p(dz)}\\
\ &\lesssim\ h^2\bigl(L+jhL_H\bigr)\sqrt{d+p}.
\end{align*}

The term proportional to $L$ enters only through the two endpoint differences, whereas the sum of corrected errors contributes $jhL_H$.
For $b>0$, on the event $\{W_j>b\}$, it holds $|e^{-W_j/p}-1|>1-e^{-b/p}$.
Markov's inequality and a union over $1\leq j\leq E$ establish
\begin{equation}\label{eq:acceptance-maximal}
 \left[\frac{\mathbb P_{\bar\mu}(M_E>b)}E\right]^{1/p}
\ \lesssim\ 
 \frac{h^2\bigl(L+EhL_H\bigr)\sqrt{d+p}}{1-e^{-b/p}}.
\end{equation}

Now fix $C,R\geq0$ and $\delta\in(0,1/2)$ as in the proposition, and suppose $Eh\leq C\Lambda_R/\sqrt m$.
Choose $b=1/16$ and $p=\lceil\log(16E/\delta)\rceil$.
Then $p\asymp\log(E/\delta)$ and $1-e^{-b/p}\gtrsim1/p$. The step-size bound \eqref{eq:pu-step-size}, with a sufficiently large implicit constant, ensures
\[
 h^2p\bigl(L+C\Lambda_RL_H/\sqrt m\bigr)\sqrt{d+p}\quad\text{is sufficiently small}.
\]
It also implies the smallness condition of \Cref{lem:acceptance-local}.
Choose this smallness so that the upper bound in \eqref{eq:acceptance-maximal}, including its universal constant, is at most $e^{-1}$. Raising to the $p$-th power proves $\mathbb P_{\hat\mu}(M_E>b)\leq Ee^{-p}\leq\delta/16$.

Set $D=\bigl\{z:\mathbb P_z(M_E>1/16)\leq1/16\bigr\}$.
The finite unadjusted path is continuous in its initial point and noises, so $D$ is measurable.
Conditioning and Markov's inequality imply $\hat\mu(D^\comp)\leq\delta$, as in the conditioning argument of \cite[Sec.~3.1]{chen2023}.
For $z\in D$, \eqref{eq:acceptance-threshold} yields
\[
 \mathbb P_z\left(\left(\bigcap\nolimits_{k=1}^E A_k\right)^\comp\right)
\ \leq\ 1-e^{-1/16}+\mathbb P_z(M_E>1/16)\ \leq\ \tfrac18.
\]
\end{proof}

\subsection{High-probability Adjusted Dobrushin Contraction over Epoch}

Now, we are prepared to prove \Cref{lem:MAKLApuDob} via combination of \Cref{prop:obabo-dobrushin,prop:pu-acceptance}, and subsequently \Cref{thm:main} via application of \Cref{thm:mixlocDobHopf}.

\begin{proof}[Proof of \Cref{lem:MAKLApuDob}]
As in \Cref{def:MAKLApu}, let $\Pi$ be the transition kernel of an epoch of $E\in\mathbb N$ steps of Metropolis-adjusted kinetic Langevin with persistent uniforms with friction $\gamma>0$ and step size $h>0$.
We first show \eqref{eq:puDobsep}, which results from a simple coupling argument.
Let $D\subseteq\R^{2d}$ and $z,\tilde z\in D$.
Couple the Gaussian random variables underlying epochs of OBABO steps started in $z$ and $\tilde z$ such that
\begin{equation}\label{eq:coup}
    Z_u\ =\ \Phi_{h,\gamma}^{\xi_E,\xi'_E}\circ\cdots\circ\Phi_{h,\gamma}^{\xi_1,\xi'_1}(z)\quad\text{and}\quad \tilde Z_u\ =\ \Phi_{h,\gamma}^{\tilde\xi_E,\tilde\xi'_E}\circ\cdots\circ\Phi_{h,\gamma}^{\tilde\xi_1,\tilde\xi'_1}(\tilde z)
\end{equation}
are a maximal coupling of $\Pi_u(z,\cdot)$ and $\Pi_u(\tilde z,\cdot)$, meaning that
\[ \TV\bigl(\Pi_u(z,\cdot),\Pi_u(\tilde z,\cdot)\bigr)\ =\ \P(Z_u\ne\tilde Z_u). \]
Such a coupling exists by first maximally coupling the endpoint distributions and then, conditionally on the coupled endpoints, sampling the Gaussian random variables from their respective regular conditional distributions.
Each resulting noise block consists of independent Gaussian random variables.
Hence, for any coupling $(U_0,\tilde U_0)$ of two $\Unif(0,1)$ random variables that is independent of the pair of noise blocks, the outputs $Z=Z_E$ and $\tilde Z=\tilde Z_E$ of the iterative procedure described in \Cref{def:MAKLApu} form a coupling of $\Pi(z,\cdot)$ and $\Pi(\tilde z,\cdot)$.
By the coupling characterization of total variation, hence,
\[ \TV\bigl(\Pi(z,\cdot),\Pi(\tilde z,\cdot)\bigr)\ \leq\ \P(Z\ne\tilde Z). \]
Using that $Z=Z_u$ on $\bigcap_{k=1}^EA_k$ and $\tilde Z=\tilde Z_u$ on $\bigcap_{k=1}^E\tilde A_k$, where $(A_k)_{1\leq k\leq E}$ and $(\tilde A_k)_{1\leq k\leq E}$ denote the respective acceptance events, the right hand side of the previous display satisfies
\begin{align*}
\P(Z\ne\tilde Z)
\ &\leq\ \P\left(Z_u\ne\tilde Z_u\,;\,\bigcap\nolimits_{1\leq k\leq E}(A_k\cap\tilde A_k)\right)\ +\ \P\left(\left(\bigcap\nolimits_{1\leq k\leq E}(A_k\cap\tilde A_k)\right)^\comp\right) \\
&\leq\ \TV\bigl(\Pi_u(z,\cdot),\Pi_u(\tilde z,\cdot)\bigr)\ +\ 2\,\sup_{z\in D}\mathbb P_z\left(\left(\bigcap\nolimits_{1\leq k\leq E}A_k\right)^\comp\right),
\end{align*}
where we used maximality of the coupling $(Z_u,\tilde Z_u)$ in the last step.

Now, suppose $\mu$ satisfies \Cref{ass:1} and let $\gamma=2\sqrt m$, $\delta\in(0,1/2)$.
For $R\geq0$, let $D_R$ and $\Lambda_R$ be as in \Cref{prop:obabo-dobrushin}.
By \Cref{prop:obabo-dobrushin,prop:pu-acceptance}, there exists
\[ h^{-1}\ \asymp\ \bigl(L^{1/2}+L_H^{1/2}m^{-1/4}\Lambda_R^{1/2}\bigr)\bigl(d+R+\log(E/\delta)\bigr)^{1/4}\bigl(\log^{1/2}(E/\delta)+\Lambda_R^{1/4}\bigr) \]
and $E\asymp\Lambda_R/(h\sqrt m)$ such that:
\begin{enumerate}[(i)]
\item For all $z,\tilde z\in D_R$,
\[ \TV\bigl(\Pi_u(z,\cdot),\Pi_u(\tilde z,\cdot)\bigr)\ \leq\ \frac12. \]

\item There exists a measurable $D'\subseteq\mathbb R^{2d}$ such that
\[ \hat\mu(D'^\comp)\ \leq\ \frac\delta2\quad\text{and}\quad\sup_{z\in D'}\mathbb P_z\left(\left(\bigcap\nolimits_{1\leq k\leq E}A_k\right)^\comp\right)\ \leq\ \frac18. \]
\end{enumerate}
Set $D=D_R\cap D'$ with a suitable $R\asymp d+\log(1/\delta)$ such that, by \cite[Lem.~5.4]{Nawaf26persistence} and (ii),
\[ \hat\mu(D^\comp)\ \leq\ \hat\mu(D_R^\comp)+\hat\mu(D'^\comp)\ \leq\ e^{-\frac12(R-3d)} + \frac\delta2\ \leq\ \delta. \]
The resulting epoch length satisfies
\[ E\ \lesssim\ \bigl(\kappa^{1/2}\Lambda+L_H^{1/2}m^{-3/4}\Lambda^{3/2}\bigr)\bigl(d+\log\delta^{-1}+\Lambda\bigr)^{1/4}\bigl(\log\delta^{-1}+\Lambda\bigr)^{1/2} \]
with $\Lambda$ as in the statement of the lemma.
Inserting (i) and (ii) into \eqref{eq:puDobsep} finishes the proof.
\end{proof}

\begin{proof}[Proof of \Cref{thm:main}]
Suppose $\mu$ satisfies \Cref{ass:1} and, as in \Cref{def:MAKLApu}, let $\Pi$ be the transition kernel of an epoch of $E\in\mathbb N$ steps of Metropolis-adjusted kinetic Langevin with persistent uniforms with friction $\gamma=2\sqrt m$ and step size $h>0$.
Let $\eps,\delta\in(0,1/2)$ and set
\[ \rho\ =\ \frac14\quad\text{and} \quad\tilde\delta\ =\ \frac{\eps\,\delta}{4\bigl\lceil\rho^{-1}\log(2/\eps)\bigr\rceil+1}. \]
By \Cref{lem:MAKLApuDob}, there exists a measurable $\tilde D\subseteq\R^{2d}$, a step size $h>0$, and an epoch length
\begin{equation}\label{eq:pfE}
    E\ \lesssim\ \bigl(\kappa^{1/2}\Lambda+L_H^{1/2}m^{-3/4}\Lambda^{3/2}\bigr)\bigl(d+\log\tilde\delta^{-1}+\Lambda\bigr)^{1/4}\bigl(\log\tilde\delta^{-1}+\Lambda\bigr)^{1/2}
\end{equation}
with $\Lambda$ as in the statement of the theorem, such that
\begin{equation}\label{eq:resdelta}
    \sup_{z,\tilde z\in\tilde D}\TV\bigl(\Pi(z,\cdot),\Pi(\tilde z,\cdot)\bigr)\ \leq\ \frac34\quad\text{and}\quad\hat\mu(\tilde D^\comp)\ \leq\ \tilde\delta.
\end{equation}
Invoking \Cref{thm:mixlocDobHopf} shows \eqref{eq:mixtimeEpoch} for a suitable measurable $D\subseteq\tilde D$ such that $\hat\mu(D^\comp)\leq\delta$ by our choice of $\tilde\delta$.
The complexity guarantee \eqref{eq:cc} follows from \eqref{eq:mixtimeEpoch} via \eqref{eq:pfE}.
\end{proof}

\subsection{Mixing from Warm Starts}

To complement the result on uniform mixing of Metropolis-adjusted kinetic Langevin in high-probability domains, we show mixing from warm start via the corresponding general framework \Cref{thm:mixlocDobwarm}.

\begin{theorem}[Mixing of Metropolis-adjusted Kinetic Langevin from Warm Starts]\label{thm:mainwarm}
Suppose $\mu$ satisfies \Cref{ass:1} and, as in \Cref{def:MAKLApu}, let $\Pi$ be the transition kernel of an epoch of $E\in\mathbb N$ steps of Metropolis-adjusted kinetic Langevin with persistent uniforms with friction $\gamma=2\sqrt m$ and step size $h>0$.
Let $\eps\in(0,1/2)$, $p\in(1,\infty]$, and $\hat\nu\in\mathcal P(\R^{2d})$ be an $L^p$-warm initial distribution.
Then, there exist suitable $E,h$ such that
\begin{equation}\label{eq:mixtimeEpochwarm}
    \tmix(\eps,\hat\nu)\ =\ \mathcal O(\log\eps^{-1}).
\end{equation}
The resulting complexity in terms of potential gradient evaluations to achieve $\eps$-accuracy in total variation from start in $\hat\nu$ is
\[ \mathcal O\left(\Lambda^{9/4}\bigl(\kappa^{1/2}+L_H^{1/2}m^{-3/4}\bigr)\Bigl(d+\tfrac{p}{p-1}\log\bigl(\beta_p(\hat\nu)/\eps\bigr)\Bigr)^{1/4}\Bigl(\tfrac{p}{p-1}\log\bigl(\beta_p(\hat\nu)/\eps\bigr)\Bigr)^{1/2}\log\eps^{-1}\right), \]
where $\Lambda=\log\max\bigl(e,\kappa,L_H,L_H/m^{3/2},d,\frac{p}{p-1}\log\bigl(\beta_p(\hat\nu)/\eps\bigr)\bigr)$.
\end{theorem}

\begin{proof}[Proof of \Cref{thm:mainwarm}]
Suppose $\mu$ satisfies \Cref{ass:1} and, as in \Cref{def:MAKLApu}, let $\Pi$ be the transition kernel of an epoch of $E\in\mathbb N$ steps of Metropolis-adjusted kinetic Langevin with persistent uniforms with friction $\gamma=2\sqrt m$ and step size $h>0$.
Let $\eps\in(0,1/2)$, $p\in(1,\infty]$, and $\hat\nu\in\mathcal P(\R^{2d})$ be $L^p$-warm.
Set
\[ \rho\ =\ \frac14\quad\text{and} \quad\tilde\delta\ =\ \frac{1}{\bigl\lceil\rho^{-1}\log(2/\eps)\bigr\rceil}\left(\frac{\eps}{4\beta_p(\hat\nu)}\right)^{\frac{p}{p-1}}. \]
By \Cref{lem:MAKLApuDob}, there exists a measurable $\tilde D\subseteq\R^{2d}$, a step size $h>0$, and an epoch length
\begin{equation}\label{eq:pfEwarm}
    E\ \lesssim\ \bigl(\kappa^{1/2}\Lambda+L_H^{1/2}m^{-3/4}\Lambda^{3/2}\bigr)\bigl(d+\log\tilde\delta^{-1}+\Lambda\bigr)^{1/4}\bigl(\log\tilde\delta^{-1}+\Lambda\bigr)^{1/2}
\end{equation}
with $\Lambda$ as in the statement of the theorem, such that
\begin{equation}\label{eq:resdeltawarm}
    \sup_{z,\tilde z\in\tilde D}\TV\bigl(\Pi(z,\cdot),\Pi(\tilde z,\cdot)\bigr)\ \leq\ \frac34\quad\text{and}\quad\hat\mu(\tilde D^\comp)\ \leq\ \tilde\delta.
\end{equation}
Invoking \Cref{thm:mixlocDobwarm} shows \eqref{eq:mixtimeEpochwarm}.
The complexity guarantee follows from \eqref{eq:mixtimeEpochwarm} via \eqref{eq:pfEwarm}, using $\log\tilde\delta^{-1}\lesssim\frac{p}{p-1}\log\bigl(\beta_p(\hat\nu)\eps^{-1}\bigr)$.
\end{proof}

\subsection*{Acknowledgements}
I thank Nawaf Bou-Rabee for sharing with me his work \cite{Nawaf26persistence}, which inspired the application of the general coupling framework to Metropolis-adjusted kinetic Langevin.
The author gratefully acknowledges funding by the Deutsche Forschungsgemeinschaft (DFG, German Research Foundation) – CRC 1720 – 539309657.

\subsection*{Statement on AI Use}

Publicly available AI tools were used for literature review and in proving \Cref{prop:obabo-dobrushin,prop:pu-acceptance}.

\printbibliography
\appendix

\section{Proofs of the Contraction Estimates}\label{app:A}

We prove the one-step estimate of \Cref{lem:readable-local} and the
two supporting estimates used in the proof of
\Cref{prop:obabo-dobrushin}. The first supporting estimate controls
energy, velocity tails and the change of $\mathcal F_h$ when leapfrog
is replaced by velocity reversal at large intermediate velocities.
The second bounds entropy after deterministic initialization.
Their proofs follow the one-step calculation below.

Throughout this section $m=1$ and $\gamma=2$. We use the normalized
parameters and the definitions of $\varphi$, $\mathcal H$,
$\mathcal C$, $G_h$, $\mathcal F_h$, $O$, $K_h$, $T_r$ and $t_r$
from \Cref{sec:obabo-contraction-proof}. The quantity $\Lambda_R$
retains its original-model value and satisfies
\eqref{eq:readable-log-controls}. All implicit constants are universal
unless a dependence is stated. Put
$\varrho=e^{-h}$ and $\sigma^2=1-e^{-2h}$.
For a probability distribution $p$ on phase space with $\mathcal H(p)<\infty$, define
\[
 \mathcal E_x(p)\ =\ \KL(p_x\mid\mu),\qquad
 \mathcal E_v(p)\ =\ \int\KL\bigl(p(dv\mid x)\mid\varphi\bigr)\,p_x(dx).
\]
The entropy chain rule reads $\mathcal H=\mathcal E_x+\mathcal E_v$.
For an input $(X,v_-)$ to leapfrog, write
\begin{equation}\label{eq:readable-physical-laws}
 a = h/2,\quad V = v_--a\nabla U(X),\quad Y = X+hV,\quad
 v_+ = V-a\nabla U(Y),\quad M = (X+Y)/2.
\end{equation}
In the full OBABO update this input has distribution $pO$.

\begin{lemma}[Energy and Velocity Bounds]\label{lem:readable-localization}
Suppose $\mu$ satisfies \Cref{ass:1} with $m=1$, and take $\gamma=2$.
Let $R\geq0$, $H(z)\leq R$, $0<\epsilon\leq1/2$, $T\geq0$, and $h>0$.
Suppose $h^2L+h^4L_H^2R_*$ is sufficiently small. Choose
\begin{equation}\label{eq:unified-exit-radius}
 R_*\ \gtrsim\ R+d+\log\frac{1+T}\epsilon,\qquad
 r\ =\ 3\sqrt{R_*},\qquad v_*\ =\ 2r,
\end{equation}
with a sufficiently large universal implicit constant.
With probability at least $1-\epsilon$, every unadjusted-chain endpoint
at time $nh\leq T$ has energy at most $R_*$, and every intermediate velocity in these
steps has $|V|\leq r$.

Define the reflected leapfrog map and comparison kernel by
\[
 \Phi_h^\#(x,v)\ =\ 
 \begin{cases}\Phi_h(x,v)&\text{if $\bigl|v-a\nabla U(x)\bigr| \leq v_*,$}\\
 S(x,v)&\text{otherwise,}\end{cases}
 \qquad\text{and}\quad K_h^\#\ =\ O\Phi_h^\#O.
\]
The chain started at $z$, taking one unadjusted step and then steps with
$K_h^\#$, differs from the unadjusted chain using the same Gaussian variables by time
$T$ with probability at most $\epsilon$.
Let $p$ be the distribution at step $k$ of the unadjusted chain when $Lh\leq1$,
or of the comparison chain. If $(k+1)h\leq T$, then
\[
 T_r(p)\ \lesssim\ M_T\sqrt\epsilon,\qquad t_r\ \lesssim\ d\sqrt\epsilon,
 \qquad M_T\ =\ e^T\bigl[R+(d+2)(T+1)\bigr].
\]
The implicit constants cover the next OU half-step and unadjusted update.
For the comparison chain and $k\geq1$ with $(k+1)h\leq T$,
\begin{equation}\label{eq:unified-transfer}
 \mathcal F_h(pK_h^\#)-\mathcal F_h(pK_h)
\ \lesssim\ \Bigl[1+M_T+d\bigl(1+\log(1/h)+T\bigr)\Bigr]\sqrt\epsilon.
\end{equation}
\end{lemma}

\begin{lemma}[Starting Entropy]\label{lem:readable-smoothing}
Let $\mu(dx)\propto e^{-U(x)}\,dx$, where $U\in C^2(\mathbb R^d)$,
$U(0)=\min U=0$, and $I_d\leq\nabla^2U(x)\leq L\,I_d$
for all $x\in\mathbb R^d$. Take $\gamma=2$ and $R\geq0$.
Suppose $h>0$ and $h^2L$ is sufficiently small.
If $Lh\leq1$ and $j=\lceil1/h\rceil$, then
\begin{equation}\label{eq:unified-block-smoothing}
 \sup_{H(z)\leq R}\mathcal H\bigl(\delta_zK_h^j\bigr)\ \lesssim\ L(R+d).
\end{equation}
If $Lh\geq1$ and \eqref{eq:readable-log-controls} holds, then
$\mathcal H(\delta_zK_h)\lesssim(R+d)\Lambda_R$ for all $z$ with $H(z)\leq R$.
\end{lemma}

\subsection{Proof of the One-step Entropy Estimate}

We refine the free-transport calculation of
\cite[App.~A]{Nawaf26persistence}, estimate the errors introduced by
the leapfrog interpolation, and combine them with the OU dissipation
as in \cite[Prop.~4.11]{Nawaf26persistence}. We use the moment estimates under $\mu$
\eqref{eq:acceptance-score}--\eqref{eq:acceptance-trace}, whose proofs
are independent of contraction.

\paragraph{Comparison with Relative Entropy.}
In this paragraph put $X=|\nabla U|^2-\mu(|\nabla U|^2)$. Then \eqref{eq:acceptance-score} implies
$\|X\|_{L^k(\mu)}\lesssim kL\sqrt d$ for integers $k\geq2$.
Since $\mu X=0$, the exponential series has no linear term. Its
remaining terms satisfy
\[
 \log\mu\bigl(e^{uX}\bigr)\ \lesssim\ u^2L^2d
 \quad\text{for all $u$ with $|u|L\sqrt d$ sufficiently small.}
\]
Entropy duality applied to $uX$ and $-uX$ implies
$|p_xX|-\mathcal E_x(p)/u\lesssim uL^2d$ for $u>0$ with $uL\sqrt d$ sufficiently small. For $\varepsilon>0$, taking
$u=h^2/(8\varepsilon)$ proves
\begin{equation}\label{eq:readable-centered-boundary}
 \bigl|p(G_h)-\mu(G_h)\bigr|-\varepsilon\mathcal H(p)
\ \lesssim\ \varepsilon^{-1}h^4L^2d,
 \quad\text{provided $h^2L\sqrt d/\varepsilon$ is sufficiently small.}
\end{equation}
The bound $\bigl|\mathcal C(p)\bigr|\leq2\sqrt{\mathcal E_x(p)\mathcal E_v(p)}$
\cite[Lem.~3.1--3.2]{Lu26entropy} implies $|\mathcal C(p)|\leq\mathcal H(p)$.
Fixing, for example,
$\varepsilon=1/4$ and taking $\eta$ sufficiently small proves \eqref{eq:unified-functional-comparison}.

\paragraph{Change of $\mathcal C$ during a Leapfrog Step.}
Let $(X,v_-)\sim p_-$ have finite entropy, and define $V,Y,v_+,M$
by \eqref{eq:readable-physical-laws}, without a preceding OU step. Put
\[
 p_+\ =\ \Law(Y,v_+),\qquad
 p_s\ =\ \Law\bigl(M+sV,V-s\nabla U(M+sV)\bigr)\quad\text{for all $-a\leq s\leq a$.}
\]
Thus $p_{-a}=p_-$, $p_a=p_+$, and $p_0=\Law(M,V)$.
The interpolating maps are compositions of shears with unit Jacobian
and at most linear growth, as are their inverses. Their changes in $H$
have at most quadratic growth. The density change-of-variables formula
and finite second moments ensure finite entropy at every stage.

For a conditional velocity distribution $\nu$ with $\KL(\nu\mid\varphi)<\infty$, let $w=\nu v$ and
$\Sigma=\nu\bigl[(v-w)(v-w)^\top\bigr]$. For a symmetric positive-definite matrix $G$, the
quadratic function
\[
 f(v)\ =\ \tfrac12(v-w)^\top(I-G)(v-w)+w\cdot v
                  -\tfrac12|w|^2+\tfrac12\log\det G
\]
satisfies $\varphi(e^f)=1$, by the Gaussian integral. Entropy duality
therefore implies
\begin{equation}\label{eq:readable-covariance}
 \Sigma:(I-G)\ \leq\ 2\KL(\nu\mid\varphi)-|w|^2-\log\det G.
\end{equation}
Here $A:C=\operatorname{tr}(A^\top C)$ denotes the Frobenius product.
Following \cite[Lem.~5.1--5.3]{Lu26entropy}, replace the test function
in \cite[App.~A.1]{Nawaf26persistence} by $f$, with $G$ the
Hessian of the convex transport potential. The coefficient on $\mathcal E_v$ is then $2$.
The term $-|w|^2$ cancels the mean-velocity term in \cite[App.~A.3]{Nawaf26persistence}.
The approximation of the optimal transport map and the singular-Hessian argument there
are unchanged. The truncation and Gaussian-convolution argument in the proof of
\cite[Lem.~4.7]{Nawaf26persistence} extends the bound to every finite-entropy input
with the positive integral coefficient reduced from $6$ to $2$.
The trapezoidal identity of \cite[Lem.~4.8]{Nawaf26persistence} adds $1$.
After combining the kicks as in \cite[Lem.~4.3, Prop.~4.9]{Nawaf26persistence}
and dropping the nonnegative terms, the coefficient in \cite[Eq.~(4.6)]{Nawaf26persistence} is $3$
in place of $7$.

To express the estimate in the interpolated velocities $V-s\nabla U(M+sV)$,
let
$\bar p_s=\Law(M+sV,V)$. Translation of a conditional velocity leaves
its differential entropy unchanged, and its Gaussian energy changes by
\[
 \mathcal E_v(\bar p_s)-\mathcal E_v(p_s)
\ =\ \E\left[sV\cdot \nabla U(M+sV)-\tfrac12s^2\bigl|\nabla U(M+sV)\bigr|^2\right].
\]
Insert this exact identity into the coefficient-$3$ bound. The added term
$3\int_{-a}^a(\cdot)\,ds-a\bigl[(\cdot)_{-a}+(\cdot)_a\bigr]$ is exactly
$\E(K_1+K_2)$, where
\begin{align*}
 K_1(q,v)\ &=\ 3\int_{-a}^a s\,v\cdot \nabla U(q+sv)\,ds
                 -a^2v\cdot\bigl[\nabla U(q+av)-\nabla U(q-av)\bigr],\\
 K_2(q,v)\ &=\ \tfrac12a^3\Bigl[\bigl|\nabla U(q+av)\bigr|^2+\bigl|\nabla U(q-av)\bigr|^2\Bigr]
                 -\tfrac32\int_{-a}^a s^2\bigl|\nabla U(q+sv)\bigr|^2\,ds.
\end{align*}
Consequently
\begin{align}\label{eq:readable-physical-corrector}
 \mathcal C(p_+)-\mathcal C(p_-)
\ &\leq\ 3\int_{-a}^a\mathcal E_v(p_s)\,ds
       -a\bigl[\mathcal E_v(p_-)+\mathcal E_v(p_+)\bigr]\notag\\
 &\quad-a\bigl[\mathcal E_x(p_-)+\mathcal E_x(p_+)\bigr]
       +\E\bigl(K_1(M,V)+K_2(M,V)\bigr).
\end{align}

\paragraph{Conditional Velocity Entropy between the Leapfrog Endpoints.}
We adapt the conditional-independence argument of
\cite[Prop.~3.1]{Nawaf26persistence} to the interpolated
leapfrog velocities. Under $h^2L\leq1/8$, we show
\begin{equation}\label{eq:unified-kinetic-bridge}
 \sup_{|s|\leq a}\mathcal E_v(p_s)
       -4\bigl[\mathcal E_v(p_-)+\mathcal E_v(p_+)\bigr]
\ \lesssim\ h^4L^2d+h^6L^3\mathcal H(p_-).
\end{equation}
Let $P=\Law(X,Y)$ and introduce the reference probability measure
\[
 d\mathsf R(x,y)\ =\ Z_h^{-1}
   e^{-\bigl[U(x)+U(y)\bigr]/2-|y-x|^2/(2h^2)}\,dx\,dy.
\]
Strong convexity ensures $0<Z_h<\infty$. The map
$(X,v_-)\mapsto(X,Y)$ has Jacobian $h^d$, and the log-density of the transformed $\hat\mu$ differs from $\log\bigl(d\mathsf R/dx\,dy\bigr)$ by a function
of at most quadratic growth. Hence $\KL(P\mid\mathsf R)<\infty$.
For $0<t<1$, let $q_t=(1-t)X+tY$ be the position at fraction $t$ of
the free-transport step. Under both $P$ and $\mathsf R$, attach a further random
variable with the same conditional distribution
\[
 W\mid X,Y\sim\mathcal N\bigl(q_t,t(1-t)h^2I_d\bigr).
\]
Completing squares expresses the reference joint density, up to its
normalizing constant, as
\[
 e^{-U(x)/2-|x-w|^2/(2th^2)}
 e^{-U(y)/2-|y-w|^2/(2(1-t)h^2)}.
\]
Thus the reference makes $X$ and $Y$ conditionally independent given
$W$. Write $P_A,\mathsf R_A$ for the marginals of any specified variable
$A$ under the two joint constructions. For each $W=w$, relative entropy
of the conditional joint distribution with respect to this product reference is
at least the sum of the two conditional marginal entropies. Integration
in $w$ and the entropy chain rule therefore imply
\[
 \KL(P_{X,W}\mid\mathsf R_{X,W})
 +\KL(P_{Y,W}\mid\mathsf R_{Y,W})
\ \leq\ \KL(P_{X,Y,W}\mid\mathsf R_{X,Y,W})
       +\KL(P_W\mid\mathsf R_W).
\]
The two constructions use the same conditional distribution of $W$ given $(X,Y)$,
so $\KL(P_{X,Y,W}\mid\mathsf R_{X,Y,W})=\KL(P\mid\mathsf R)$.
Forgetting $W$ cannot increase either entropy on the left. Hence
\[
 \KL(P_X\mid\mathsf R_X)+\KL(P_Y\mid\mathsf R_Y)
\ \leq\ \KL(P\mid\mathsf R)+\KL(P_W\mid\mathsf R_W)
\ \leq\ \KL(P\mid\mathsf R)+\KL(P_{q_t}\mid\mathsf R_{q_t}).
\]
The last inequality is data processing through the common Gaussian
noise. Set
\[
 D(t)\ =\ \KL(P\mid\mathsf R)-\KL(P_{q_t}\mid\mathsf R_{q_t}).
\]
The preceding inequality reads $D(t)\leq D(0)+D(1)$. By the chain
rule, $D(t)$ is the mean conditional entropy of the interpolated velocity
given $q_t$, relative to its conditional distribution under $\mathsf R$.

At position $q=q_t$, the interpolated velocity is
$w=V+h(1/2-t)\nabla U(q)$. Let $\rho_{t,q}$ be its conditional distribution under
$\mathsf R$. Differentiation of the pair density shows that
\begin{equation}\label{eq:readable-reference-score}
 \begin{aligned}
 \nabla_w\log\frac{d\rho_{t,q}}{d\varphi}(w)
\ &=\ \tfrac h2\Bigl(t\bigl[\nabla U(X)-\nabla U(q)\bigr]-(1-t)\bigl[\nabla U(Y)-\nabla U(q)\bigr]\Bigr),\\
 \left|\nabla_w\log\frac{d\rho_{t,q}}{d\varphi}\right|
\ &\leq\ h^2L|w|+h^3L\bigl|\nabla U(q)\bigr|.
 \end{aligned}
\end{equation}
Here $X=q-thV$ and $Y=q+(1-t)hV$, with
$V=w-h(1/2-t)\nabla U(q)$. These formulas also define the reference kernels at
$t=0,1$.

To compare this reference with $\varphi$, suppose $f=\log\bigl(d\rho/d\varphi\bigr)$ is
$C^1$, $\bigl|\nabla f(v)\bigr|\leq\alpha|v|+\beta$, where $0\leq\alpha\leq1/8$ and $\beta\geq0$. Set
$\psi(u)=\log\varphi\bigl(e^{uf}\bigr)$ and
$d\rho_u=e^{uf-\psi(u)}d\varphi$. The bound on $\nabla f$ implies
quadratic growth of $f$, so all the following integrals and derivatives
exist for all $|u|\leq2$, with domination from $|u|\leq3$. Applying
the Gaussian logarithmic Sobolev inequality \cite{Gross76} to $e^{uf/2}$ gives
\[
 D_u\ :=\ \KL(\rho_u\mid\varphi)\ =\ u\psi'(u)-\psi(u)
\ \leq\ u^2\bigl[\alpha^2\rho_u|v|^2+\beta^2\bigr].
\]
Entropy duality against $|v|^2/4$, whose Gaussian exponential moment is
$2^{d/2}$, bounds $\rho_u|v|^2$ by $4D_u+2d\log2$.
Absorbing $4u^2\alpha^2D_u$ leaves
$D_u\lesssim u^2\bigl(\alpha^2d+\beta^2\bigr)$.
Integration of $\bigl(\psi(u)/u\bigr)'=D_u/u^2$ on either side of zero shows that
$0\leq\psi(u)-u\psi'(0)\lesssim u^2\bigl(\alpha^2d+\beta^2\bigr)$. Since $\psi(1)=0$,
$\psi'(0)\leq0$. Define the order-two R\'enyi divergence by
$D_2(\rho\mid\lambda)=\log\int\bigl(d\rho/d\lambda\bigr)^2d\lambda$
when $\rho\ll\lambda$, and by $+\infty$ otherwise. We have shown that
\[
 D_2(\rho\mid\varphi)+D_2(\varphi\mid\rho)
\ =\ \psi(2)+\psi(-1)\ \lesssim\ \alpha^2d+\beta^2.
\]
If $\nu,\lambda,\rho$ are probability distributions with
$\KL(\nu\mid\lambda)<\infty$ and $D_2(\lambda\mid\rho)<\infty$,
entropy duality and Cauchy--Schwarz give, for every bounded measurable $f$,
\[
 \nu(f)\ \leq\ 2\KL(\nu\mid\lambda)+2\log\lambda(e^{f/2})
\ \leq\ 2\KL(\nu\mid\lambda)+D_2(\lambda\mid\rho)+\log\rho(e^f).
\]
Taking the supremum over $f$ in the variational formula for relative entropy yields
\[
 \KL(\nu\mid\rho)\ \leq\ 2\KL(\nu\mid\lambda)+D_2(\lambda\mid\rho).
\]
For each $q$, use this inequality with $(\rho,\lambda)$ equal first
to $(\rho_{t,q},\varphi)$ and then to $(\varphi,\rho_{t,q})$.
The two-sided estimate just proved and
\eqref{eq:readable-reference-score} apply with
$\alpha=h^2L$ and $\beta=h^3L\bigl|\nabla U(q)\bigr|$. Write
$G(t)=\mathcal E_v(p_{h(t-1/2)})$. Averaging over $q_t$, we obtain
\[
 \max\bigl(D(t)-2G(t),G(t)-2D(t)\bigr)\ \lesssim\ h^4L^2d+h^6L^2\E\bigl|\nabla U(q_t)\bigr|^2,
\]
by changing the reference in each direction.
Together with $D(t)\leq D(0)+D(1)$, these inequalities imply
$G(t)-4\bigl[G(0)+G(1)\bigr]\lesssim h^4L^2d+h^6L^2\sup_{0\leq u\leq1}\E\bigl|\nabla U(q_u)\bigr|^2$.
The endpoint cases are immediate.

Since $q_t=X+thV$, $|\nabla U|^2\leq2LU$, and
$V=v_--(h/2)\nabla U(X)$, the pointwise bound
$\bigl|\nabla U(q_t)\bigr|^2\leq6L H(X,v_-)$ holds when $h^2L\leq1/8$.
Entropy duality using \eqref{eq:readable-energy-mgf} at $u=1/2$ implies
$p_-(H)\leq2\mathcal H(p_-)+4d$. Hence
$\E\bigl|\nabla U(q_t)\bigr|^2\lesssim L\bigl[\mathcal H(p_-)+d\bigr]$.
Substitution now shows that
$G(t)-4\bigl[G(0)+G(1)\bigr]\lesssim h^4L^2d+h^6L^3\bigl[\mathcal H(p_-)+d\bigr]$.
The condition $h^2L\leq1/8$ absorbs the last dimensional term into $h^4L^2d$, proving
\eqref{eq:unified-kinetic-bridge}.

\paragraph{Full Relative Entropy at the Leapfrog Midpoint.}
We also need the full midpoint entropy bound
\begin{equation}\label{eq:unified-full-bridge}
 \mathcal H(p_0)-\mathcal H(p_-)-\mathcal H(p_+)\ \lesssim\ h^4(L^2+L_H^2)d.
\end{equation}
We compare the interpolated map with exact Hamiltonian flow, then pair
the resulting displacement with its negative to cancel the linear
energy and log-determinant terms.

The argument requires both $h\sqrt L$ and
$h^4\bigl(L^2+L_H^2/L\bigr)$ to be sufficiently small. Both follow from
\eqref{eq:unified-local-cap}, since $L\geq1$ and $d\geq1$, after
decreasing its universal smallness threshold. For this paragraph only,
rescale position and elapsed time by $\sqrt L$, replacing the potential
by $x\mapsto U(x/\sqrt L)$. Put $b=h\sqrt L$ and $\beta=L_H/L^{3/2}$,
and retain $U,H,\hat\mu$ for the transformed objects within this paragraph.
The upper Hessian bound is now one and the Frobenius Hessian Lipschitz
constant is $\beta$. Energy is unchanged at corresponding points.
For an initial point $z=(x,v)$ and $0\leq u\leq b$, put
$V_b=v-(b/2)\nabla U(x)$ and
\[
 q_u\ =\ x+uV_b,\qquad
 w_u\ =\ V_b+(b/2-u)\nabla U(q_u),\qquad
 \mathcal M_u(x,v)\ =\ (q_u,w_u).
\]
Let
$\Psi_u(x,v)=(x_u,v_u)$ be the exact flow solving
$\dot x_u=v_u$, $\dot v_u=-\nabla U(x_u)$ from the same
point. Subtracting the
integral equations, we obtain
\begin{align*}
 q_u-x_u\ &=\ \tfrac12u(u-b)\nabla U(x)
       +\int_0^u(u-t)\bigl[\nabla U(x_t)-\nabla U(x)\bigr]\,dt,\\
 w_u-v_u\ &=\ (b/2-u)\bigl[\nabla U(q_u)-\nabla U(x)\bigr]
       +\int_0^u\bigl[\nabla U(x_t)-\nabla U(x)\bigr]\,dt.
\end{align*}
Conservation of $H$ under exact flow implies
$|x_u-x|+|q_u-x|\lesssim b\sqrt{H(x,v)}$.
Since $|\nabla U|^2\leq2U$ and
$\|\nabla^2U\|_{\mathrm{op}}\leq1$, insertion into these
identities proves
$|\mathcal M_u-\Psi_u|\lesssim b^2\sqrt{H}$. To estimate derivatives, differentiate the
same identities, and write each difference of differentiated forces in
the form
\[
 \bigl[\nabla^2U(q)-\nabla^2U(x)\bigr]Dq
       +\nabla^2U(x)(Dq-Dx).
\]
The Frobenius norm of the first term is at most
$\beta|q-x|\|Dq\|_{\mathrm{op}}$. For the second, the variational equations
of exact flow have bounded coefficients. Gronwall's inequality implies
$\|D\Psi_u\|_{\mathrm{op}}\leq e^{|u|}$ and
$\|D\Psi_u-I\|_F\lesssim |u|\sqrt d$. Directly from the formula for $q_u$,
we also have
$\|Dq_u\|_{\mathrm{op}}\lesssim1$ and $\|Dq_u-Dx\|_F\lesssim u\sqrt d$.
Thus each differentiated force difference above has Frobenius norm
$\lesssim b(\beta\sqrt{H}+\sqrt d)$.
The coefficients and time integrals supply another factor $b$.
The explicit term $u(u-b)\nabla U(x)/2$ has derivative of Frobenius norm $\lesssim b^2\sqrt d$.
For $T_u=\mathcal M_u\circ\Psi_{-u}=\mathrm{id}+e_u$, compose the derivative
error with $D\Psi_{-u}$ to obtain
\begin{equation}\label{eq:readable-flow-errors}
 |e_u|^2\ \lesssim\ b^4H,\qquad
 \|De_u\|_F^2\ \lesssim\ b^4\bigl(d+\beta^2H\bigr),\qquad
 \|De_u\|_{\mathrm{op}}\ \lesssim\ b.
\end{equation}
For the spectral norm bound compare each derivative directly with the
identity, using only the bounded Hessian. No derivative of the Hessian
is required in these integral identities.

The maps $\mathcal M_u$, $\Psi_u$, and $T_u$ preserve volume. For small
$b$, $T_u^-:=\mathrm{id}-e_u$ is a global $C^1$ diffeomorphism: for each
$y$, the equation $z=y+e_u(z)$ is a contraction. Since
$|\nabla H|^2\leq2H$, the function $\sqrt{H}$ is $1/\sqrt2$-Lipschitz.
Together with $|e_u(z)|\lesssim b^2\sqrt{H(z)}$, this permits absorption of
the small factor and proves
$H\bigl((T_u^-)^{-1}y\bigr)\lesssim H(y)$.
Taylor's theorem for energy and the expansion of
$\log\bigl|\det(I\pm De_u)\bigr|$ imply
\begin{align*}
 &H(T_uz)+H(T_u^-z)-2H(z)\\
 &\quad-\log\bigl|\det DT_u(z)\bigr|-\log\bigl|\det DT_u^-(z)\bigr|
 \ \lesssim\ b^4\bigl[d+(1+\beta^2)H(z)\bigr].
\end{align*}
The linear energy terms and linear traces cancel between the two signs.
The remaining terms are $\lesssim |e_u|^2+\|De_u\|_F^2$.
For the determinant step one uses, for any real matrix $A$ with
$\|A\|_{\mathrm{op}}\leq1/2$,
$-\log\bigl|\det(I+A)\bigr|-\log\bigl|\det(I-A)\bigr|\lesssim\|A\|_F^2$.
No symmetry of $A=De_u$ is needed.

By volume preservation,
\[
 \int\left(\frac{d\bigl((\mathcal M_u)_\#\hat\mu\bigr)}
                   {d\hat\mu}\right)^2d\hat\mu
\ =\ \hat\mu
       \exp\bigl[H(\mathcal M_uz)-H(z)\bigr]
\ =\ \hat\mu
       \exp\bigl[H(T_uz)-H(z)\bigr].
\]
The second equality follows by substituting $z=\Psi_{-u}y$ and using
invariance of $\hat\mu$ under $\Psi_u$. Inserting the preceding
energy and determinant inequality produces the Jacobian of $T_u^-$.
Changing variables by this map and using
$H\bigl((T_u^-)^{-1}y\bigr)\lesssim H(y)$, we obtain
\[
 \log\hat\mu\exp\bigl[H(T_uz)-H(z)\bigr]
 \ \lesssim\ b^4d+\log\hat\mu\bigl(e^{\lambda H}\bigr),
\]
with $\lambda\asymp b^4(1+\beta^2)$,
where the universal comparison constant is sufficiently large.
The step-size condition can be chosen so that $\lambda\leq1/2$.
The convexity calculation in \cite[Lem.~5.4]{Nawaf26persistence} proves, more generally,
\begin{equation}\label{eq:readable-energy-mgf}
 \hat\mu\bigl(e^{uH}\bigr)\ \leq\ (1-u)^{-3d/2}\quad\text{for all $0\leq u<1$.}
\end{equation}
Indeed convexity and $U(0)=0$ imply
$\int e^{-(1-u)U}\leq(1-u)^{-d}\int e^{-U}$.
The Gaussian velocity factor contributes $(1-u)^{-d/2}$.
The same exponential-moment bound holds before rescaling.
It follows, when $b^4(1+\beta^2)$ is sufficiently small, that
$D_2\bigl((\mathcal M_u)_\#\hat\mu\mid\hat\mu\bigr)
\lesssim b^4(1+\beta^2)d$.
KL invariance under $\mathcal M_u$ and the reference-change inequality
imply, with $p_-$ and $\mathcal H$ expressed in these coordinates,
\[
 \mathcal H\bigl((\mathcal M_u)_\#p_-\bigr)-2\mathcal H(p_-)
\ \lesssim\ b^4(1+\beta^2)d.
\]
Velocity reversal establishes the same midpoint bound from $p_+$.
Average the two bounds and return to the $m=1$ units, where
$b^4(1+\beta^2)=h^4\bigl(L^2+L_H^2/L\bigr)\leq h^4(L^2+L_H^2)$. This proves
\eqref{eq:unified-full-bridge}.

\paragraph{Expectations of the Energy and Quadrature Errors.}
We first estimate Gaussian moments of the errors, then use a transport
inequality to bound their expectations under the midpoint distribution. Define
the energy remainder and midpoint entropies by
\[
 \begin{gathered}
 \mathcal R_{\mathrm{mid}}(q,v)\ =\ U(q+av)-U(q-av)-av\cdot\bigl[\nabla U(q+av)+\nabla U(q-av)\bigr],\\
 E_m\ =\ \mathcal E_v(p_0),\qquad X_m\ =\ \mathcal E_x(p_0).
 \end{gathered}
\]
We prove the two bounds used in \eqref{eq:readable-physical-corrector}:
\begin{align}
 \bigl|\E\mathcal R_{\mathrm{mid}}(M,V)\bigr|\ &\lesssim\ L_Hh^3\bigl(\sqrt{dE_m}+rE_m\bigr)+Lh^2T_r,
 \label{eq:unified-signed-remainder}\\
 \E(K_1+K_2)(M,V)-\tfrac h8\mathcal H(p_0)\ &\lesssim\ h^5(L^2+L_H^2)d
                 +L^2h^3T_r+Lh^3t_r.
 \label{eq:unified-quadrature}
\end{align}
Here $T_r=\E\bigl[|V|^2\mathbf1_{\{|V|>r\}}\bigr]$ under the arbitrary
finite-entropy midpoint distribution. It equals $T_r(p)$ when $p_-=pO$.

In this paragraph, write $J=\nabla^2U$ for the Hessian.
For a fixed position $q$ and scalar $s\in\mathbb R$, let
$A_s(v)=\bigl[J(q+sv)-J(q)\bigr]v$. We claim that
\begin{equation}\label{eq:readable-hessian-increment}
 \varphi|A_s|^2\ \lesssim\ L_H^2s^2d.
\end{equation}
For $s\ne0$, symmetry of the almost-everywhere weak third
derivatives implies
\[
 \sum_i\bigl|DJ(q+sv)[e_i]v\bigr|^2
\ =\ \bigl\|DJ(q+sv)[v]\bigr\|_F^2\ \leq\ L_H^2|v|^2,
\]
where $(e_i)$ is an orthonormal basis. Differentiating
$A_s(v)=\bigl[J(q+sv)-J(q)\bigr]v$ consequently proves
$\bigl\|DA_s(v)\bigr\|_F\lesssim |s|L_H|v|$. Gaussian integration by parts identifies
$\varphi A_s=s\varphi\nabla\Delta U(q+sv)$, whose norm is at most
$|s|L_H\sqrt d$ because $|\nabla\Delta U|\leq\sqrt d\,L_H$. Componentwise Gaussian
Poincar\'e \cite[Ch.~4]{bakry_book} bounds the variance. Both the variance
and squared mean are $\lesssim L_H^2s^2d$.
This proves \eqref{eq:readable-hessian-increment}. The weak-derivative
calculation can equivalently be justified by mollifying $U$ with
unchanged Hessian bounds and passing by dominated convergence.

The force Taylor remainder will be denoted by
$R_s(v)=\nabla U(q+sv)-\nabla U(q)-sJ(q)v$. By the fundamental theorem of calculus,
$R_s(v)=s\int_0^1 A_{us}(v)\,du$. Jensen's inequality then implies
$\varphi|R_s|^2\lesssim L_H^2s^4d$.
Direct differentiation of the scalar energy remainder shows that
\[
 \nabla_v\mathcal R_{\mathrm{mid}}(q,v)\ =\ -a^2\bigl[J(q+av)-J(q-av)\bigr]v.
\]
For $K_1$, put
\[
 D_s(v)\ =\ R_s(v)-R_{-s}(v)+s\bigl[A_s(v)+A_{-s}(v)\bigr].
\]
Differentiating its defining integral and cancelling the constant-Hessian
terms using $3\int_0^1u^2du=1$ proves
\[
 \nabla_vK_1(q,v)\ =\ a^2\left[3\int_0^1uD_{au}(v)\,du-D_a(v)\right].
\]
Thus \eqref{eq:readable-hessian-increment} and Jensen imply
\[
 \varphi\bigl|\nabla_v\mathcal R_{\mathrm{mid}}\bigr|^2\ \lesssim\ L_H^2h^6d,\qquad
 \varphi|\nabla_vK_1|^2\ \lesssim\ L_H^2h^8d.
\]
By the same formulas and Hessian Lipschitzness, the gradient Lipschitz
constants on $|v|\leq2r$ are $\lesssim L_Hh^3r$ and $\lesssim L_Hh^4r$, respectively.
The pointwise bounds needed below are
\[
 \bigl|\mathcal R_{\mathrm{mid}}\bigr|\ \lesssim\ L_Hh^3|v|^3\wedge Lh^2|v|^2,\qquad
 |K_1|\ \lesssim\ L_Hh^4|v|^3,\qquad |K_1|\ \leq\ Lh^3|v|^2.
\]

Let $\chi_r$ be an even $C^{1,1}$ radial cutoff, equal to one on
$|v|\leq r$ and zero on $|v|\geq2r$, with
$|\nabla\chi_r|\lesssim 1/r$ and
$\operatorname{Lip}(\nabla\chi_r)\lesssim 1/r^2$.
The products $\chi_r\mathcal R_{\mathrm{mid}}$ and $\chi_rK_1$ have the same orders of
Gaussian gradient energy and global gradient Lipschitz constants as
above. Indeed the product rule adds $f\nabla\chi_r$ to the gradient,
and the three second-derivative terms are bounded by the local gradient
Lipschitz bound, $|\nabla f|/r$, and $|f|/r^2$ on the annulus.
For the gradient energy, the cutoff derivative requires a bound on
$\varphi(f^2)/r^2$, with $f=\mathcal R_{\mathrm{mid}}$ or $K_1$. Gaussian Poincar\'e bounds
the centered part by $\varphi|\nabla f|^2/r^2$, and
$\varphi\mathcal R_{\mathrm{mid}}=0$ because $\mathcal R_{\mathrm{mid}}$ is odd in velocity. For $K_1$ the mean term
also satisfies $|\varphi K_1|^2/r^2\lesssim L_H^2h^8d$, because
$|\varphi K_1|\lesssim L_Hh^4d$ by \eqref{eq:readable-K1-mean} below, and $r^2\geq d$.

For any scalar $f$ with square-integrable Gaussian gradient and globally
Lipschitz gradient, and any distribution $\nu$ with $\KL(\nu\mid\varphi)<\infty$, choose an optimal quadratic-cost coupling of
$Z\sim\varphi$ and $V\sim\nu$. Taylor expansion at $Z$, followed by
Cauchy--Schwarz for its linear term, shows that
\[
 \bigl|\nu f-\varphi f\bigr|
\ \leq\ \sqrt{\varphi|\nabla f|^2}\,W_2(\nu,\varphi)
       +\tfrac12\operatorname{Lip}(\nabla f)W_2(\nu,\varphi)^2.
\]
The Gaussian transport inequality \cite{Talagrand96} therefore implies
\begin{equation}\label{eq:unified-gaussian-taylor}
 \bigl|\nu f-\varphi f\bigr|\ \leq\ 
 \sqrt{2\varphi\bigl(|\nabla f|^2\bigr)\KL(\nu\mid\varphi)}
       +\operatorname{Lip}(\nabla f)\KL(\nu\mid\varphi).
\end{equation}
Apply this conditionally given $M$ to the cutoff of $\mathcal R_{\mathrm{mid}}$, whose
Gaussian mean is zero. Conditional Jensen bounds the integral of
$\sqrt{\KL}$ by $\sqrt{E_m}$. Removing the cutoff costs
$\lesssim Lh^2T_r$, and proves \eqref{eq:unified-signed-remainder}.

For the nonzero Gaussian mean of $K_1$, set
$\bar K_1(q)=\varphi K_1(q,\cdot)$ and
$P_sf(q)=\int f(q+sv)\,\varphi(dv)=e^{(s^2/2)\Delta}f(q)$. Gaussian integration by parts identifies this mean as
\begin{equation}\label{eq:readable-K1-mean}
 \bar K_1\ =\ 3\int_{-a}^a s^2P_s\Delta U\,ds
                 -a^3(P_a+P_{-a})\Delta U.
\end{equation}
The coefficient of the constant function vanishes. Since
$\operatorname{Lip}(\Delta U)\leq L_H\sqrt d$, it follows that
$|\bar K_1|\lesssim L_Hh^4d$ and
$\operatorname{Lip}(\bar K_1)\lesssim L_Hh^3\sqrt d$.
To average this mean under $\mu$, let
$\mu_t=\mu*\mathcal N(0,2tI_d)$, the distribution of $X+\sqrt{2t}Z$ for
independent $X\sim\mu$ and $Z\sim\varphi$. For bounded Lipschitz $f$,
Gaussian differentiation and integration by parts express the derivative as
\[
 \frac d{dt}\mu_t f\ =\ \operatorname{Cov}_{\mu_t}(f,b_t),\qquad
 b_t(y)\ =\ \E\Bigl[\bigl|\nabla U(X)\bigr|^2-\Delta U(X)\mid X+\sqrt{2t}Z=y\Bigr].
\]
The estimate \eqref{eq:acceptance-trace}, conditional Jensen,
and its mean-zero identity imply $\|b_t\|_{L^2(\mu_t)}\lesssim L\sqrt d$.
For $t>0$, the conditional density of $X$ given
$X+\sqrt{2t}Z=y$ has negative log-density with Hessian at least $\bigl[1+(2t)^{-1}\bigr]I_d$.
The Poincar\'e inequality under this curvature bound \cite{Bakry85}
therefore bounds its covariance by $2t(1+2t)^{-1}I_d$.
Twice differentiating the convolution density expresses the Hessian of
its negative logarithm as
$(2t)^{-1}I_d-(4t^2)^{-1}\operatorname{Cov}(X\mid y)$.
It is consequently at least $(1+2t)^{-1}I_d$. A second application of the
same criterion bounds the Poincar\'e constant of $\mu_t$ by $1+2t$.
At $t=0$ the corresponding bound is that of $\mu$. In particular,
$\operatorname{Var}_{\mu_t}(f)\leq(1+2t)\operatorname{Lip}(f)^2$.
Cauchy--Schwarz in the covariance formula bounds the derivative by
$L\sqrt d\,\operatorname{Lip}(f)\sqrt{1+2t}$ up to a universal factor.
Integration over $0\leq t\leq s^2/2\leq1/2$ proves
\[
 \bigl|\mu P_sf-\mu f\bigr|\ \lesssim\ s^2L\sqrt d\,\operatorname{Lip}(f).
\]
Apply this to $f=\Delta U$ in \eqref{eq:readable-K1-mean}:
$|\mu \bar K_1|\lesssim L_HLh^5d$.

Applying
\eqref{eq:unified-gaussian-taylor} conditionally to $\chi_rK_1$ bounds
the difference from that mean. The transport inequality
\cite{Otto00}, using $I_d\leq\nabla^2U$, states that
$W_2\bigl((p_0)_x,\mu\bigr)^2\leq2X_m$. Together with the Lipschitz bound on
$\bar K_1$, this bounds the difference between $\E \bar K_1(M)$ and $\mu \bar K_1$. The two discarded tails cost at most $Lh^3(T_r+t_r)$, by the
explicit quadratic bound above. Altogether,
\begin{equation}\label{eq:readable-K1-bound}
 \begin{aligned}
 \bigl|\E K_1(M,V)\bigr|\ \lesssim\ {}&L_Hh^3\sqrt{dX_m}
       +L_Hh^4\bigl(\sqrt{dE_m}+rE_m\bigr)+L_HLh^5d\\
       &+Lh^3(T_r+t_r).
 \end{aligned}
\end{equation}

For $K_2$ we need a second-moment bound on the vector remainder $R_s$.
The derivative and quadratic growth bounds imply
$\operatorname{Lip}(\chi_rR_s)\lesssim L_Hs^2r$.
Applying $(a+b)^2\leq2a^2+2b^2$ under a Gaussian optimal coupling gives
\[
 \E\bigl|R_s(V)\bigr|^2\ \lesssim\ L_H^2s^4(d+r^2E_m)+L^2s^2T_r.
\]
Here $q=M$ in $R_s$, the conditional estimate is integrated over $M$, and the last term
bounds the part outside the cutoff using $\bigl|R_s(v)\bigr|\leq2L|s||v|$.
For the even difference
$F_s=\bigl|\nabla U(q+sv)\bigr|^2+\bigl|\nabla U(q-sv)\bigr|^2-2\bigl|\nabla U(q)\bigr|^2$, use the expansion
$\nabla U(q\pm sv)=\nabla U(q)\pm sJ(q)v+R_{\pm s}(v)$. For $0<|s|\leq1$,
Young's inequality implies
\[
 |F_s|\ \leq\ 2s^2\bigl|\nabla U(q)\bigr|^2+4s^2L^2|v|^2
               +3s^{-2}\bigl(\bigl|R_s(v)\bigr|^2+\bigl|R_{-s}(v)\bigr|^2\bigr).
\]
The constant part cancels from $K_2$:
$K_2=(a^3/2)\bigl[F_a-3\int_0^1u^2F_{au}du\bigr]$.
Entropy duality using \eqref{eq:readable-energy-mgf} implies
$\E\bigl|\nabla U(M)\bigr|^2\lesssim L(d+X_m)$ and
$\E|V|^2\lesssim (d+E_m)$. Integration of the preceding inequalities proves
\[
 \bigl|\E K_2(M,V)\bigr|\ \lesssim\ h^5\bigl(L(d+X_m)+L^2(d+E_m)
                       +L_H^2(d+r^2E_m)\bigr)+L^2h^3T_r.
\]
Young's inequality bounds each square-root term in
\eqref{eq:readable-K1-bound} by any prescribed small multiple of $h$ times its entropy plus a constant
multiple of $h^5(L^2+L_H^2)d$, using $h\leq1$ for the second square-root term. Also $2L_HL\leq (L^2+L_H^2)$. The coefficients of the remaining $X_m$ and $E_m$ terms, including
$L_Hh^4rE_m$, can be made arbitrarily small relative to $h$ by choosing
the universal threshold in \eqref{eq:unified-local-cap} small enough. Since
$X_m+E_m=\mathcal H(p_0)$, the total entropy contribution is at most
$h\mathcal H(p_0)/8$. The integrated $K_2$ tail
is $\lesssim L^2h^3T_r$. Adding the $K_1$ tails and using $L\geq1$
bounds their sum by a universal multiple of $L^2h^3T_r+Lh^3t_r$. This proves
\eqref{eq:unified-quadrature}.

\paragraph{Combining Leapfrog with the OU Half-steps.}
Substitution of \eqref{eq:unified-kinetic-bridge} into
\eqref{eq:readable-physical-corrector} first implies
\begin{align*}
 &\mathcal C(p_+)-\mathcal C(p_-)-\tfrac{23}2h\bigl[\mathcal E_v(p_-)+\mathcal E_v(p_+)\bigr]\\
 &\quad+\tfrac h2\bigl[\mathcal E_x(p_-)+\mathcal E_x(p_+)\bigr]-\E(K_1+K_2)
 \ \lesssim\ h^5L^2d+h^7L^3\mathcal H(p_-).
\end{align*}
Now use \eqref{eq:unified-quadrature} and
\eqref{eq:unified-full-bridge}, and take $h^2L$ sufficiently small to
absorb the last entropy term into $h\mathcal H(p_-)/8$. The resulting bound is
\begin{equation}\label{eq:readable-Verlet-stage}
 \begin{aligned}
 \mathcal C(p_+)-\mathcal C(p_-)&-13h\bigl[\mathcal E_v(p_-)+\mathcal E_v(p_+)\bigr]+\tfrac h4\mathcal E_x(p_-)\\ &\lesssim\ h^5(L^2+L_H^2)d+L^2h^3T_r+Lh^3t_r.
 \end{aligned}
\end{equation}

Apply the OU calculation of \cite[Prop.~4.11]{Nawaf26persistence}
with $a_1=13$, $a_2=1/4$, retaining the centered gradient-square term.
Let
$p_-=pO$, $p_+=p_-\Phi_h$, and $p'=p_+O$, and write
$\delta=1-\varrho^2$,
$V_0=\mathcal E_v(p)$,
$V_1=\mathcal E_v(p_+)$. The unchanged OU identities
\cite[Lem.~4.2]{Nawaf26persistence} are
$\mathcal E_x(pO)=\mathcal E_x(p)$, $\mathcal E_v(pO)\leq\varrho^2\mathcal E_v(p)$, and $\mathcal C(pO)=\varrho\mathcal C(p)$.
Together with the leapfrog estimate and volume preservation, they imply
\begin{align*}
 &\mathcal H(p')\ \leq\ \mathcal H(p)-\delta(V_0+V_1)
                     +p'(G_h)-p(G_h)+\E\mathcal R_{\mathrm{mid}}(M,V),\\
 &\mathcal C(p')-\varrho^2\mathcal C(p)-13h(V_0+V_1)+\tfrac14\varrho h\mathcal E_x(p)\ \lesssim\ h^5(L^2+L_H^2)d+L^2h^3T_r+Lh^3t_r.
\end{align*}
In the first line we used the exact energy identity
$\Delta H_h(X,v_-)=G_h(Y)-G_h(X)+\mathcal R_{\mathrm{mid}}(M,V)$. OU half-steps do not move positions.
Thus its boundary difference cancels exactly in $\mathcal F_h$.

For $h\leq1/4$, we have $h\leq\delta\leq2h$ and $\varrho\geq3/4$.
We may therefore fix the universal $\eta>0$ small enough that, uniformly in $h$,
\[
 13\eta h\ \leq\ \delta/8,\qquad
 2\eta\delta\ \leq\ \varrho h/16.
\]
For $\E\mathcal R_{\mathrm{mid}}(M,V)$, Young's inequality shows that the square-root
contribution in \eqref{eq:unified-signed-remainder} minus
$\varepsilon_0hE_m$ is $\lesssim h^5L_H^2d$, with
a fixed universal $\varepsilon_0>0$ sufficiently small relative to $\eta$.
Taking $L_Hh^2r$ sufficiently small absorbs the linear $E_m$
term into $\varepsilon_0hE_m$. Then \eqref{eq:unified-kinetic-bridge}, $\mathcal H(pO)\leq
\mathcal H(p)$, and small $h^2L$ imply
\[
 \E\mathcal R_{\mathrm{mid}}(M,V)-\tfrac\delta4(V_0+V_1)
       -\tfrac{\eta\varrho h}{8}\mathcal E_x(p)\ \lesssim\ h^5(L^2+L_H^2)d+Lh^2T_r.
\]
More explicitly, before comparing coefficients with the OU decrease,
the entropy terms are bounded by $8\varepsilon_0h(V_0+V_1)+\varepsilon_0h\bigl[\mathcal E_x(p)+V_0\bigr]$.
The bounds $\delta\geq h$ and $\varrho\geq3/4$, together with the choice of $\varepsilon_0$, imply the displayed bound.

Finally, the bound on $\mathcal C(p)$ implies
$-\eta\delta\mathcal C(p)\leq\delta V_0/2+2\eta^2\delta\mathcal E_x(p)$.
Substituting these bounds into the change of $\mathcal F_h$ leaves dissipation at least
$\delta V_0/8+\eta\varrho h\mathcal E_x(p)/16$.
After fixing $\eta$ and $\varepsilon_0$, choose a universal $c>0$ small enough that this is at least
$2ch\mathcal H(p)$. All remaining implicit constants are then universal. Therefore
\[
 \mathcal F_h(p')-\mathcal F_h(p)+2ch\mathcal H(p)
\ \lesssim\ h^5(L^2+L_H^2)d+\bigl(Lh^2+L^2h^3\bigr)T_r+Lh^3t_r.
\]
Use \eqref{eq:unified-functional-comparison} to
bound $ch\mathcal F_h(p)-2ch\mathcal H(p)\lesssim h^5L^2d$.
This proves \eqref{eq:unified-local-recursion} and completes
\Cref{lem:readable-local}.

\subsection{Energy Control and Comparison with a Reflected Chain}

We prove \Cref{lem:readable-localization} by stopping the exponential
Lyapunov estimate of \cite[Lem.~5.5]{Nawaf26persistence}, replacing
the stepwise union bound in \cite[Lem.~5.6]{Nawaf26persistence}. For the comparison chain,
the energy moments and densities are bounded as in \cite[Lem.~5.7 and~5.11]{Nawaf26persistence}.
We also compare $\mathcal C$ on the common part of the two distributions.

\paragraph{Probability of Large Energy or Intermediate Velocity.}
Let $(Z_k)_{k\geq0}$ be the unadjusted chain with $Z_0=z$. Set
$\mathcal W=H-G_h$ and $\mathcal V=\mathcal W+x\cdot v/5$.
Since $L\geq1$, decreasing the universal threshold on $h^2L$ ensures
$h\leq1/8$. The comparison and drift estimates
\cite[Eq.~(5.1), Lem.~5.5]{Nawaf26persistence} then ensure
$H/2\leq\mathcal V\leq2H$.
For a step started at $z$ with $H(z)\leq R_*$, let $\xi$ denote the
standard Gaussian used by its first OU half-step and set
$A=\bigl\{\sigma|\xi|\leq\sqrt{R_*}\bigr\}$. On $A$,
\[
 |V|\ \leq\ \sqrt{2R_*}+h\sqrt{LR_*/2}+\sqrt{R_*}\ <\ 3\sqrt{R_*}.
\]
The smallness of $h^4L_H^2R_*$ therefore bounds the remainder
$(L_Hh^3/12)|V|^3$ by a sufficiently small universal multiple of
$h|V|^2$. This is the condition used in \cite{Nawaf26persistence}.
Its Gaussian integration establishes, for fixed universal $\theta,c,C>0$,
\[
 \E\bigl[e^{\theta\mathcal V(Z_1)}\mathbf1_A\bigr]
\ \leq\ \exp\Bigl(\theta\bigl[(1-ch)\mathcal V(z)+Chd\bigr]\Bigr).
\]
If $\mathcal V(z)\geq Cd/c$, the exponent does not increase.
Otherwise $\mathcal V(z)\lesssim d$, so $e^u-1\leq ue^u$ for $u\geq0$ implies
\[
 \E\bigl[e^{\theta\mathcal V(Z_1)}\mathbf1_A\bigr]
       -e^{\theta\mathcal V(z)}\ \lesssim\ he^{Cd},
\]
where the universal constant $C$ may increase. In particular, the
coefficient of the current exponential moment remains one.

Stop $e^{\theta\mathcal V(Z_k)}$ at the first energy exit or failed noise event.
If the energy exits before any noise failure, retain the exit value thereafter.
Otherwise set the stopped process to zero. Its conditional expected
increment is $\lesssim he^{Cd}$, and its initial value is at most $e^{2\theta R}$.
At an energy exit before any failed noise event, its value is at least
$e^{\theta R_*/2}$. By Markov's inequality, the energy-exit probability by time $nh\leq T$
is $\lesssim e^{-cR_*}\bigl(e^{CR}+Te^{Cd}\bigr)$.
For failed noise events, the Gaussian exponential moment and
$\sigma^2\leq2h$ bound the total probability by
\[
 \frac Th\exp\bigl(Cd-cR_*/h\bigr)
\ \lesssim\ T\exp\bigl(Cd-cR_*\bigr).
\]
Here we used $R_*\geq1$ and $h\leq1/8$ to absorb $h^{-1}$ into
half of the negative exponential, adjusting the universal constants.
The sufficiently large constant in \eqref{eq:unified-exit-radius}
makes the sum at most $\epsilon$, uniformly in $h$.
Outside this event all endpoint energies are at most $R_*$ and all
intermediate velocities satisfy $|V|\leq r$.

Couple the unadjusted and comparison chains using the same noises.
Since $r<v_*$, induction shows that they agree at every stage outside
the exceptional event. Thus, for either chain, each event $|V|>r$ has
probability at most $\epsilon$, including the next unadjusted step ending by time $T$.

\paragraph{Energy Moments and Density Bounds for the Two Chains.}
To convert these probabilities into weighted tails by Cauchy--Schwarz,
we first bound the energy second moments of both chains.

The comparison map $\Phi_h^\#$ preserves volume
\cite[Lem.~5.7]{Nawaf26persistence}: insert at the midpoint the
involution that reverses velocity on $\bigl\{|v|>v_*\bigr\}$, and compose
with the position and velocity shears on either side.
The same energy calculation establishes
\[
 \mathcal W\circ\Phi_h^\#\ \leq\ (1+h)\mathcal W,
 \qquad \mathcal W\circ\Phi_h\ \leq\ (1+Lh^2)\mathcal W.
\]
Indeed $\mathcal W\circ\Phi_h-\mathcal W=\mathcal R_{\mathrm{mid}}(M,V)$.
On an accepted comparison step, $|V|\leq v_*=2r$ and
$\bigl|\mathcal R_{\mathrm{mid}}\bigr|\leq(L_Hh^3v_*/12)|V|^2$.
On a reflected step $\mathcal W$ is unchanged. For an unrestricted step,
$\bigl|\mathcal R_{\mathrm{mid}}\bigr|\leq Lh^2|V|^2/8$.
Using $U\geq|\nabla U|^2/(2L)$ and completing the square, we obtain
\[
 \mathcal W\ \geq\ \tfrac12\bigl(1-h^2L/4\bigr)|V|^2,
 \qquad \bigl(1-h^2L/4\bigr)H\ \leq\ \mathcal W\ \leq\ H.
\]
Since $L_Hh^2v_*\lesssim\sqrt{h^4L_H^2R_*}$, these bounds justify
both recurrence coefficients after reducing the universal threshold.
They also compare $\mathcal W$ and $H$ with universal constants.
For any current distribution $p$ with finite energy second moment, an OU
half-step satisfies
\[
 \|\mathcal W\|_{L^2(pO)}
\ \leq\ \|\mathcal W\|_{L^2(p)}+h(d+2).
\]
Indeed, expansion of the Gaussian quadratic velocity energy shows that
$O(\mathcal W^2)(x,v)\leq\bigl[\mathcal W(x,v)+\sigma^2(d+2)/2\bigr]^2$.
Minkowski's inequality and $\sigma^2\leq2h$ imply the displayed bound,
as in \cite[Lem.~5.11]{Nawaf26persistence}. The position part
of $\mathcal W$ is unchanged. Consequently, a full comparison step satisfies
\[
 \|\mathcal W\|_{L^2(pK_h^\#)}
\ \leq\ (1+h)\|\mathcal W\|_{L^2(p)}+3h(d+2).
\]
Iterating with $(1+h)^n\leq e^{nh}$ proves
$\|H\|_{L^2(p)}\lesssim M_T$ for every distribution of the comparison chain within
the time horizon. The single initial unadjusted step changes only the
universal constant, since $h^2L$ is small. If $Lh\leq1$, the unadjusted
factor $1+Lh^2$ is at most $1+h$, so the same bound holds for
distributions of the unadjusted chain. An additional OU half-step and unadjusted update
change only the universal constant. Thus the two outputs in
\eqref{eq:unified-transfer} satisfy
\[
 \max\bigl(\|H\|_{L^2(pK_h)},\|H\|_{L^2(pK_h^\#)}\bigr)\ \lesssim\ M_T,
 \qquad \bigl(\E|V|^4\bigr)^{1/2}\ \lesssim\ M_T.
\]

Cauchy--Schwarz with the event from the preceding paragraph now proves
$T_r\lesssim M_T\sqrt\epsilon$ for both chains in their stated regimes.
Since the chosen radius has $r^2$ sufficiently large compared with $d$,
the Gaussian exponential moment also implies
$\varphi\bigl(|v|^2\mathbf1_{\{|v|>r\}}\bigr)\lesssim d e^{-cr^2}$.
Increasing the radius constant makes this at most $d\sqrt\epsilon$.

The density after the first unadjusted step is at most
$(2\pi h\sigma^2)^{-d}$, by the Gaussian-coordinate computation
\cite[Lem.~4.13]{Nawaf26persistence}. A volume-preserving deterministic
map leaves this upper bound unchanged. An OU half-step increases an
upper Lebesgue-density bound by at most $\varrho^{-d}=e^{dh}$, because
convolution does not increase a density supremum. Thus every comparison
output under consideration has density at most
\begin{equation}\label{eq:readable-density-cap}
 D_*\ =\ (2\pi h\sigma^2)^{-d}e^{2dT}.
\end{equation}
Since $\sigma^2\geq h$,
$\log D_*\lesssim d\bigl[1+\log(1/h)+T\bigr]$.
Together with the energy moments, this also ensures finite entropy at
all distributions of the comparison chain and their unadjusted updates.

\paragraph{Comparison of $\mathcal F_h$ for Distributions Coupled to Agree.}
We extend the entropy mixture estimate of
\cite[Lem.~5.8]{Nawaf26persistence} to $\mathcal F_h$, using
the log-Jacobian argument \cite[Sec.~2 and~4]{McCann97}
to compare the optimal transport maps on the common measure.

Let $P,Q$ have
finite relative entropy with respect to $\hat\mu$, $\|H\|_{L^2(P)},\|H\|_{L^2(Q)}\leq M$, and let $Q$
have Lebesgue density at most $D$. Suppose they share a subprobability
measure $\zeta\leq P,Q$ of mass $q=1-\alpha$, with $0\leq\alpha\leq1/2$.
The case $\alpha=0$ is immediate. For $0<\alpha\leq1/2$, the cited
mixture calculation proves
\[
 \mathcal H(Q)-\mathcal H(P)
\ \leq\ \alpha\mathcal H\bigl((Q-\zeta)/\alpha\bigr)
                   -\alpha\log\alpha-q\log q.
\]
Indeed, convexity bounds the $Q$ mixture from above, while the chain
rule bounds the $P$ mixture from below up to its Bernoulli entropy.
The density of $(Q-\zeta)/\alpha$ is at most $D/\alpha$.
Moreover, strong convexity bounds $\log\int e^{-H}$ by $d\log(2\pi)$,
while Cauchy--Schwarz bounds $\int H\,d\bigl(Q-\zeta\bigr)$ by $M\sqrt\alpha$.
Substituting these bounds and using $-\alpha\log\alpha-q\log q\lesssim \sqrt\alpha$ proves, with
$\log_+D=\max(\log D,0)$,
\begin{equation}\label{eq:readable-entropy-transfer}
 \mathcal H(Q)-\mathcal H(P)
\ \lesssim\ \bigl[1+M+d+\log_+D\bigr]\sqrt\alpha.
\end{equation}

To compare $\mathcal C(P)$ and $\mathcal C(Q)$, let $T_0,T_1$ be the optimal transport maps from $P_x,Q_x$ to
$\mu$, respectively, and put $\beta=\zeta_x/q$.
Since $\beta\leq P_x/q,Q_x/q$, its images
$\nu_i=(T_i)_\#\beta$ satisfy $\nu_i\leq\mu/q$, hence
$\KL(\nu_i\mid\mu)\leq\log(1/q)$.
The log-Jacobian proof of displacement convexity
\cite[Thm.~4.4 and proof of Thm.~2.2]{McCann97}, applied to these two gradient maps
\emph{from the same source}, implies
\[
 \KL\left(\left(\tfrac{T_0+T_1}{2}\right)_\#\beta\,\Big|\,\mu\right)
\ \leq\ \tfrac12\KL(\nu_0\mid\mu)+\tfrac12\KL(\nu_1\mid\mu)
                            -\tfrac18\beta\bigl|T_0-T_1\bigr|^2.
\]
In the change-of-variables formulas, the source-density terms cancel
pointwise before integration, so finite entropy of $\beta$ is not needed.
The quadratic term follows from strong convexity of $U$, and the
Hessian inequality from concavity of $\log\det$. Absolute
continuity of $\beta$ suffices: restrict both optimal transport maps to a common
set of full $\beta$-measure on which they are injective with positive Alexandrov
Hessians. Strict monotonicity of an injective convex gradient makes
their average injective there, with positive Hessian, so the same
change-of-variables formula applies.
Nonnegativity of the midpoint KL therefore implies
\begin{equation}\label{eq:readable-common-map}
 \int\bigl|T_0-T_1\bigr|^2\,d\zeta_x
\ \lesssim\ q\log(1/q)\ \leq\ \alpha.
\end{equation}

On the common measure, the difference contributing to $\mathcal C(Q)-\mathcal C(P)$ is
$\int v\cdot\bigl(T_0(x)-T_1(x)\bigr)\,d\zeta$ and is
$\lesssim\sqrt{M\alpha}$ using \eqref{eq:readable-common-map}.
For the residual measures, Cauchy--Schwarz and the energy bound imply
\[
 \int\bigl(|x|^2+|v|^2\bigr)\,d(P-\zeta)\ \leq\ 2M\sqrt\alpha,
 \qquad \int|T_0(x)|^2\,d(P-\zeta)\ \leq\ \sqrt{\alpha\mu|x|^4},
\]
with the same estimates for $Q,T_1$. The second inequality uses that
$(T_0)_\#(P_x-\zeta_x)$ is a submeasure of $\mu$ with mass $\alpha$. Strong convexity and
\eqref{eq:readable-energy-mgf} imply
$\mu|x|^4\lesssim (d+1)^2$: use $|x|^2\leq2U$ and apply the bound
$\mu(e^{uU})\leq(1-u)^{-d}$ at $u=1/(d+1)$, together with
$u^2U^2\leq4e^{uU}$.
Thus the integrals defining $\mathcal C$ over $P-\zeta$ and $Q-\zeta$ are
$\lesssim(M+d+1)\sqrt\alpha$. Finally $G_h\leq H$ under $h^2L\leq4$, so the
difference of its expectations is $\lesssim M\sqrt\alpha$.
Together with \eqref{eq:readable-entropy-transfer}, this proves
\[
 \mathcal F_h(Q)-\mathcal F_h(P)
\ \lesssim\ \bigl[1+M+d+\log_+D\bigr]\sqrt\alpha.
\]

Apply this with $P=pK_h$, $Q=pK_h^\#$ and the common output measure
on $\{|V|\leq v_*\}$. The preceding estimates verify
$\alpha\leq\epsilon$, $M\lesssim M_T$ and $D=D_*$, proving
\eqref{eq:unified-transfer} and \Cref{lem:readable-localization}.

\subsection{Entropy after Deterministic Initialization}

For \Cref{lem:readable-smoothing}, one step suffices when $Lh\geq1$
\cite[Lem.~4.13]{Nawaf26persistence}. For $Lh\leq1$, we avoid its
$\log(1/h)$ term by comparing consecutive steps with the same integrator
with zero force,
using the Gaussian coordinates of \cite[Eq.~(17)--(18)]{Bussi07}.

\paragraph{Steps of Total Time between One and Two.}
Let $K_h^0$ denote OBABO with zero force and the same OU parameters as
$K_h$. For an actual transition from $z=(x,v)$, the position and conditional
velocity are
\[
 Y\ =\ x+h\varrho v-\tfrac12h^2\nabla U(x)+h\sigma\xi,\qquad
 v'\ =\ \varrho\bigl((Y-x)/h-(h/2)\nabla U(Y)\bigr)+\sigma\xi',
\]
where $\xi,\xi'$ are independent standard Gaussians. The free
transition omits the two force terms without changing the conditional
covariances. The chain rule and the Gaussian relative-entropy formula identify
\begin{equation}\label{eq:readable-one-step-reference}
 \KL\bigl(K_h(z,\cdot)\mid K_h^0(z,\cdot)\bigr)
\ =\ \frac{h^2}{8\sigma^2}
       \Bigl[\bigl|\nabla U(x)\bigr|^2+\varrho^2\E_z\bigl|\nabla U(Y)\bigr|^2\Bigr].
\end{equation}
For $H(z)\leq R$ and $Lh\leq1$, the unadjusted-chain moment estimate
proved above bounds the expected energies by a universal multiple of $R+d$ up to time two,
including the position $Y$ in each step. As $|\nabla U|^2\leq2LU$ and
$\sigma^2\geq h$, \eqref{eq:readable-one-step-reference} contributes
$\lesssim hL(R+d)$ per transition to the relative entropy of the two
path distributions. Summing by the chain rule, then applying data processing to
the endpoint, proves
\[
 \KL\bigl(\delta_zK_h^j\mid\delta_z(K_h^0)^j\bigr)\ \lesssim\ L(R+d)
 \quad\text{for all $1\leq jh\leq2$.}
\]
For one position--velocity coordinate pair, the free transition
matrix and covariance are
\[
 A\ =\ \begin{pmatrix}1&h\varrho\\0&\varrho^2\end{pmatrix},\qquad
 C_h\ =\ \sigma^2\begin{pmatrix}h^2&h\varrho\\h\varrho&1+\varrho^2\end{pmatrix}.
\]
The corresponding endpoint covariance is
$C_{j,h}=\sum_{k=0}^{j-1}A^kC_h(A^k)^\top$. Summing these geometric
series, we find, for $t=jh$,
\[
 \det C_{j,h}
\ =\ th\coth h\bigl(1-e^{-4t}\bigr)
       -\bigl(h/\sinh h\bigr)^2\bigl(1-e^{-2t}\bigr)^2
\ \gtrsim\ 1\quad\text{for all $1\leq t\leq2$.}
\]
Indeed $h\coth h\geq1$ and $h/\sinh h\leq1$ imply
$\det C_{j,h}\geq t\bigl(1-e^{-4t}\bigr)-\bigl(1-e^{-2t}\bigr)^2
\geq2e^{-2t}\bigl(1-e^{-2t}\bigr)$ for $1\leq t\leq2$.
The $d$ coordinate pairs are independent, so the logarithm of the free endpoint
density is $\lesssim d$, uniformly in $h$ and the starting point.

Let $P=\delta_zK_h^j$ and $N=\delta_z(K_h^0)^j$.
The distribution $P$ has finite second moments and $N$ is a nondegenerate Gaussian,
so its quadratic log-density is integrable under $P$.
The relative-entropy chain rule for a change of reference reads
\[
 \KL(P\mid\hat\mu)-\KL(P\mid N)-P H
\ =\ P\log\frac{dN}{dz}+\log\int e^{-H}
\ \lesssim\ d.
\]
Combining this identity with the path-entropy and energy estimates proves
\eqref{eq:unified-block-smoothing}. The displayed upper bounds also
justify finiteness of the entropies in this computation.

\paragraph{One Step When $\mathbf{Lh\geq1}$.}
By \cite[Lem.~4.13]{Nawaf26persistence},
\[
 \mathcal H(\delta_zK_h)
\ \leq\ \E_z H(Z_1)+d\log\frac1{2e h^2(1-h)}.
\]
When $h^2L$ is sufficiently small, we have $h\leq1/8$, and the same
elementary energy calculation bounds $\E_z H(Z_1)$ by a universal multiple of $R+d$.
If $Lh\geq1$, then
$\log(1/h)\leq\log L\leq\Lambda_R$. This proves \Cref{lem:readable-smoothing}.

\subsection{Choice of Epoch Parameters}\label{sec:obabo-epoch-parameters}

We verify the choices and the iteration in the proof of
\Cref{prop:obabo-dobrushin}. Let $c>0$ be the decay constant in
\Cref{lem:readable-local}. Choose a universal $K\geq1$ large enough
to dominate the constants in the three lemmas and the geometric sum.
With $\bar h,E,T$ as in \eqref{eq:readable-epoch-choice}, set
$\epsilon=e^{-K(T+\Lambda_R)}$ and $R_*=C_A(R+d)\Lambda_R$.
The inequalities \eqref{eq:readable-log-controls} imply
\[
 R+d+\log\frac{1+T}\epsilon\ \lesssim_A\ (R+d)\Lambda_R.
\]
Thus a sufficiently large $C_A$ ensures \eqref{eq:unified-exit-radius}.
Next choose $c_A$ small. For $h\leq\bar h$,
\[
 h^4\bigl[(L^2+L_H^2)d+L_H^2R_*\bigr]\ \lesssim\ C_Ac_A^4,
\]
since $R+d,\Lambda_R\geq1$ and
$\bigl(\sqrt L+\sqrt{L_H}\bigr)^4\geq (L^2+L_H^2)$.
With $r^2=9R_*$, this implies all step-size conditions in the three
lemmas. The same choice makes $h\leq1$, so $A\geq2$ ensures
$E\geq\lceil1/h\rceil$ and $(E+2)h\leq T$.
The initial elapsed time is $t=jh$ with $1\leq t\leq2$ when $Lh\leq1$,
and $t=h\leq1$ otherwise. Consequently $Nh\geq Eh-2$ in both cases.

For the unadjusted sequence, \Cref{lem:readable-smoothing} implies
$\mathcal H(p_0)\lesssim L(R+d)$. For the comparison sequence it implies
$\mathcal H(p_0)\lesssim(R+d)\Lambda_R$.
The tail bounds inserted into \eqref{eq:unified-local-recursion}
produce an error $\lesssim h\bigl[h^4(L^2+L_H^2)d+P\sqrt\epsilon\bigr]$ with
\[
 \mathcal H(p_0)\ \lesssim\ L(R+d)\Lambda_R,\qquad
 P\ =\ L^2e^T(R+d)(T+\Lambda_R+1).
\]
Indeed $M_T\lesssim e^T(R+d)(T+1)$, and the remaining coefficients
are bounded by this $P$ using $h\leq1$ and $L\geq1$.
For the comparison sequence, dividing \eqref{eq:unified-transfer} by
$h$ costs at most $L$, and $\log(1/h)\leq\log L\leq\Lambda_R$.
Its error is therefore covered by the same $hP\sqrt\epsilon$.
All implicit constants here are universal.

Iterating \eqref{eq:unified-local-recursion} with these error bounds
and using $\sum_{k\geq0}(1-ch)^k=(ch)^{-1}$ bounds
$\mathcal F_h(p_N)-e^{-cNh}\max(\mathcal F_h(p_0),0)$ by a universal
multiple of $h^4(L^2+L_H^2)d+P\sqrt\epsilon$. The positive part is needed because
$\mathcal F_h(p_0)$ can be negative. The two sides of
\eqref{eq:unified-functional-comparison} imply
\[
 \mathcal H(p_N)-4e^{-cNh}\mathcal H(p_0)
\ \lesssim\ h^4(L^2+L_H^2)d+P\sqrt\epsilon.
\]
By \eqref{eq:readable-log-controls},
$\log\bigl(1+\mathcal H(p_0)\bigr)\lesssim\Lambda_R$ and
$\log(1+P)\lesssim T+\Lambda_R$.
The coupling in \Cref{lem:readable-localization}
and $Nh\geq Eh-2$ therefore establish
\[
 \begin{aligned}
 \TV\bigl(\delta_zK_h^E,p_N\bigr)\ &\leq\ \epsilon,\\
 \KL(p_N\mid\hat\mu)-4e^{-c(Eh-2)}\mathcal H(p_0)
\ &\lesssim\ \Bigl(h\bigl(\sqrt L+\sqrt{L_H}\bigr)\bigl((R+d)\Lambda_R\bigr)^{1/4}\Bigr)^4+P\sqrt\epsilon.
 \end{aligned}
\]
First choose $K$ large enough to make $P\sqrt\epsilon$ sufficiently
small. Next choose $A$ large enough to make the transient term small,
using $Eh\geq A\Lambda_R$ and $\Lambda_R\geq1$. Then choose $C_A$ for
the radius and $c_A$ for all step-size conditions and the
fourth-order error. These successive universal choices ensure the
bounds $\epsilon\leq1/16$ and $\KL(p_N\mid\hat\mu)\leq1/32$ used
in the main proof.

\section{Proofs of the Moment Estimates}\label{app:B}

We prove \Cref{lem:acceptance-local} with the notation and original, unscaled parameters of \Cref{sec:obabo-acceptance-proof}.
The argument has four parts. Centered gradient and Laplacian moments control the logarithmic density ratio under Gaussian convolution. This transfers Laplacian moments to shifted positions and bounds the density in midpoint coordinates. Gaussian estimates for the corrected remainder then complete the weighted moment bound.
The gradient and Laplacian estimates also enter \Cref{app:A}, but their proofs are independent of contraction.

\subsection{Gradient and Laplacian Moments}

\begin{lemma}[Centered Gradient and Laplacian Moments]\label{lem:readable-gibbs}
Let $\mu(dx)\propto e^{-U(x)}\,dx$, where $U\in C^2(\R^d)$ satisfies $m\,I_d\leq\nabla^2U(x)\leq L\,I_d$ for $0<m\leq L$ and all $x\in\R^d$.
With $F,X$ as above, set $J=\Delta U-F$.
For every integer $q\geq1$,
\begin{alignat}{2}
 \mu(F)\ &\leq\ Ld,
 &\qquad &\|X\|_{L^q(\mu)}
 \ \lesssim\ L\bigl(\sqrt{dq}+q\bigr)\ \lesssim\ qL\sqrt d,
 \label{eq:acceptance-score}\\
 \mu(J)\ &=\ 0,
 &\qquad &\|J\|_{L^q(\mu)}
       +\bigl\|\Delta U-\mu(\Delta U)\bigr\|_{L^q(\mu)}
 \ \lesssim\ qL\sqrt d.
 \label{eq:acceptance-trace}
\end{alignat}
Moreover, $\mu\bigl(e^{2s|X|}\bigr)\lesssim1$ if $sL\sqrt d$ is sufficiently small.
\end{lemma}

\begin{proof}
We first obtain exponential moments of the centered gradient norm by transport and exponential tilting. We convert these to centered moments of $F$, then use the change of variables $x\mapsto x+t\nabla U(x)$ to estimate $J$ and the Laplacian.
For $\KL(\nu\mid\mu)<\infty$, strong convexity and entropy duality imply finite second moments of $\nu$.
Let $T_\#\mu=\nu$ be the optimal transport map for squared Euclidean cost \cite{Brenier91}.
Here $DT$ denotes the Alexandrov Hessian of its convex potential.
The weak Jacobian formula \cite[Thm.~4.4]{McCann97} and the integration-by-parts argument of \cite[App.~A.1]{Nawaf26persistence}, with source $\mu$ and displacement $T-x$, yield
\[
 \KL(\nu\mid\mu)\ =\ \mu\bigl[U\circ T-U-\log\det DT\bigr]\quad\text{and}\quad
 \mu\operatorname{tr}(DT-I_d)\ \leq\ \mu\bigl[\nabla U\cdot(T-x)\bigr].
\]
Finite second moments and finite relative entropy make the energy difference
and logarithmic density ratio integrable, hence also $\log\det DT$.
In the second inequality, the singular part of the distributional Hessian of the transport potential has nonnegative trace.
Radial cutoffs have vanishing boundary terms by the second moments.
Thus $\log\det A\leq\operatorname{tr}(A-I_d)$ yields
\begin{equation}\label{eq:readable-bregman}
 \mu\bigl[U\circ T-U-\nabla U\cdot(T-x)\bigr]\ \leq\ \KL(\nu\mid\mu).
\end{equation}
The convexity and smoothness inequality
\[
 \bigl|\nabla U(y)-\nabla U(x)\bigr|^2
\ \leq\ 2L\bigl[U(y)-U(x)-\nabla U(x)\cdot(y-x)\bigr]
\]
follows by taking a gradient step of length $1/L$ for $U(\cdot)-\nabla U(x)\cdot(\cdot)$ and comparing with its minimum at $x$.
Therefore, for $Y=|\nabla U|$, $\bigl|\nu Y-\mu Y\bigr|\leq\sqrt{2L\,\KL(\nu\mid\mu)}$.
For every $t\in\R$, the exponential tilt $d\nu_t\propto e^{t(Y-\mu Y)}d\mu$ has finite entropy because $Y$ grows linearly and $\mu$ has Gaussian tails.
Its entropy identity reads
\[
 \log\mu\bigl(e^{t(Y-\mu Y)}\bigr)
\ =\ t\,\nu_t\bigl(Y-\mu Y\bigr)-\KL(\nu_t\mid\mu)
\ \leq\ \sup\nolimits_{k\geq0}\bigl(|t|\sqrt{2Lk}-k\bigr)\ =\ Lt^2/2.
\]
Set $Y_0=Y-\mu Y$. Chernoff's inequality with $t=\pm r/L$ implies $\mu\bigl(|Y_0|>r\bigr)\leq2e^{-r^2/(2L)}$ for $r>0$. Cavalieri's principle then expresses the moment as a tail integral:
\begin{align*}
 \|Y_0\|_{L^q(\mu)}
\ &=\ \left(q\int_0^\infty r^{q-1}\mu\bigl(|Y_0|>r\bigr)\,dr\right)^{1/q}\\
\ &\leq\ \left(2q\int_0^\infty r^{q-1}e^{-r^2/(2L)}\,dr\right)^{1/q}
\ \lesssim\ \sqrt{Lq}.
\end{align*}
Integration by parts shows that $\mu F=\mu(\Delta U)\leq Ld$.
Since
\[
 \mu Y\ \leq\ \sqrt{Ld},\qquad X\ =\ 2(\mu Y)Y_0+Y_0^2-\mu(Y_0^2),
\]
the triangle inequality and the preceding bound at orders $q$ and $2q$ prove \eqref{eq:acceptance-score}.
For $t\geq0$, since $n^n/n!\leq e^n$, monotone convergence also permits the expansion
\[
 \mu\bigl(e^{t|X|}\bigr)\ =\ 1+\sum_{n\geq1}\frac{t^n\mu\bigl(|X|^n\bigr)}{n!}\ \lesssim\ 1
 \quad\text{when }tL\sqrt d\text{ is sufficiently small}.
\]

For $J=\Delta U-F$, use the diffeomorphism $T_t(x)=x+t\nabla U(x)$, $|t|L\leq1/2$.
Its derivative is uniformly positive definite. For each $y$, the minimizer of $|x|^2/2+tU(x)-y\cdot x$ is the unique solution of $T_t(x)=y$.
The change-of-variables formula is $1=\mu\Bigl(e^{U-U\circ T_t}\det\bigl(I_d+t\nabla^2U\bigr)\Bigr)$.
The bounds
\[
 U(T_tx)-U(x)\ \leq\ tF(x)+Lt^2F(x)/2,
 \qquad t\Delta U-\log\det\bigl(I_d+t\nabla^2U\bigr)\ \lesssim\ t^2L^2d
\]
therefore imply
\begin{equation}\label{eq:readable-compensated-J}
 \log\mu\bigl(e^{tJ-Lt^2F/2}\bigr)\ \lesssim\ t^2L^2d.
\end{equation}
To remove the term in $F$, use the recurrence of \cite[Lem.~9, App.~B.1]{chen2023}, with $f=U$ and $\Upsilon=Ld$:
\[\begin{aligned}
 \mu(F^{k+1})
\ &=\ \mu(F^k\Delta U)+2k\mu\bigl(F^{k-1}\nabla U^\top\nabla^2U\nabla U\bigr)\\
\ &\leq\ L(d+2k)\mu(F^k)\quad\text{for all $k\geq1$.}
\end{aligned}\]
These integrations by parts use only $C^2$ regularity, and Gaussian tails control their cutoff terms. Starting with $\mu F\leq Ld$, iteration and expansion of the exponential show that $\mu\bigl(e^{uF}\bigr)\leq(1-2Lu)^{-d/2}\lesssim1$ when $u\geq0$ and $uLd$ is sufficiently small.
By Cauchy--Schwarz and \eqref{eq:readable-compensated-J} at $2t$,
\[
 \mu\bigl(e^{tJ}\bigr)\ \leq\ 
 \mu\bigl(e^{2tJ-2Lt^2F}\bigr)^{1/2}\mu\bigl(e^{2Lt^2F}\bigr)^{1/2}\ \lesssim\ 1
 \quad\text{if }|t|L\sqrt d\text{ is sufficiently small}.
\]
The bounds for both signs imply $\mu\bigl(e^{t|J|}\bigr)\lesssim1$ for sufficiently small $tL\sqrt d>0$. Choose such a $t\asymp(L\sqrt d)^{-1}$ and integrate Chernoff's tail bound to obtain $\|J\|_{L^q(\mu)}\lesssim qL\sqrt d$.
Finally, $\mu J=0$ and $\Delta U-\mu(\Delta U)=J+X$ prove \eqref{eq:acceptance-trace}.
\end{proof}

\subsection{Gaussian Shifts and the Midpoint Density}
To estimate the density after the midpoint change of variables, we need Laplacian moments under Gaussian shifts of $\mu$.
The following $L^2(\mu)$ bound on the density ratio permits this transfer by H\"older's inequality.

\begin{lemma}[Gaussian Convolution]\label{lem:readable-heat}
Suppose $\mu$ satisfies the assumptions of \Cref{lem:readable-gibbs}.
Let $T\geq0$, $\rho=d\mu/dx$, $\rho_t=e^{t\Delta}\rho$ and $\ell_t=\log(\rho_t/\rho)$.
If $TL\sqrt d$ is sufficiently small, then
\[
 \mu\bigl(|\ell_t|>b\bigr)\ \leq\ e^{-3b}\quad\text{for all $0\leq t\leq T$ and $b\geq1$.}
\]
In particular, if $h^2L\sqrt d$ is sufficiently small,
\begin{equation}\label{eq:acceptance-heat}
 \mu\bigl(e^{2|\ell_t|}\bigr)\ \lesssim\ 1,\qquad
 \|\rho_t/\rho\|_{L^2(\mu)}\ \lesssim\ 1\quad\text{for all $0\leq t\leq s$.}
\end{equation}
\end{lemma}

\begin{proof}
We bound moments of the derivative of $\ell_t$ under $\rho_t$, then integrate up to the first crossing of a fixed level by $|\ell_t|$.
The case $T=0$ is immediate.
For $t>0$, put $U_t=-\log\rho_t$.
In the convolution formula
\[
 \rho_t(y)\ =\ \frac1{Z(4\pi t)^{d/2}}
 \int e^{-U(x)-|y-x|^2/(4t)}\,dx,\qquad Z\ =\ \int e^{-U},
\]
the posterior $\pi_{t,y}$ has potential $V(x)=U(x)+|y-x|^2/(4t)$ with curvature bounds $\alpha=m+(2t)^{-1}$ and $\beta=L+(2t)^{-1}$.
Talagrand's inequality \cite[Lem.~4.1(S2)]{Nawaf26persistence}, with lower curvature $\alpha$, and exponential tilting as in the preceding proof bound the centered logarithmic moment-generating function of $u\cdot x$ by $\theta^2|u|^2/(2\alpha)$. Differentiating twice at $\theta=0$ proves $\operatorname{Var}_{\pi_{t,y}}(u\cdot x)\leq|u|^2/\alpha$.
Conversely, integration by parts and Cauchy--Schwarz imply
\[
 |u|^4
\ =\ \Bigl[\mathbb E_{\pi_{t,y}}
       \bigl(u\cdot(x-\mathbb E_{\pi_{t,y}}x)\bigr)\bigl(u\cdot\nabla V\bigr)\Bigr]^2
\ \leq\ \operatorname{Var}_{\pi_{t,y}}(u\cdot x)\,\beta|u|^2,
\]
using $\mathbb E\nabla V\nabla V^\top=\mathbb E\nabla^2V$.
The posterior has Gaussian tails, so these integrations are justified.
Differentiating the convolution normalizer thus yields
\[
 \nabla^2U_t(y)\ =\ \frac{I_d}{2t}-\frac{\operatorname{Cov}_{\pi_{t,y}}(x)}{4t^2},
 \qquad
 \frac m{1+2mt}\,I_d\ \leq\ \nabla^2U_t\ \leq\ \frac L{1+2Lt}\,I_d\ \leq\ L\,I_d.
\]
Thus \Cref{lem:readable-gibbs} applies to $\rho_t\,dx$, with the same upper Hessian bound $L$.
For $J_t=\Delta U_t-|\nabla U_t|^2$, its conclusion is
\begin{equation}\label{eq:readable-heat-J}
 \left(\int|J_t|^q\rho_t\,dx\right)^{1/q}\ \lesssim\ qL\sqrt d.
\end{equation}

The heat equation reads $\partial_t\ell_t=-J_t$ and $\ell_0=0$.
Since $\rho$, $\nabla\rho=-\rho\nabla U$ and $\nabla^2\rho=\rho\bigl(\nabla U\nabla U^\top-\nabla^2U\bigr)$ are bounded and continuous, Gaussian convolution extends these identities continuously to $t=0$ at each point.
If $\bigl|\ell_T(x)\bigr|>b$, integration up to its first crossing of $b$ in absolute value shows that
\[
 b\ \leq\ \int_0^T\mathbf1_{\{|\ell_t(x)|\leq b\}}\bigl|J_t(x)\bigr|\,dt.
\]
On $\bigl\{|\ell_t|\leq b\bigr\}$, $\rho\leq e^b\rho_t$.
Jensen, Tonelli and \eqref{eq:readable-heat-J} therefore imply
\begin{align*}
 \mu\bigl(|\ell_T|>b\bigr)^{1/q}
\ &\leq\ \left[\frac{T^{q-1}}{b^q}\int_0^T\!\int
       \mathbf1_{\{|\ell_t|\leq b\}}|J_t|^q\rho\,dx\,dt\right]^{1/q}\\
\ &\lesssim\ e^{b/q}qTL\sqrt d/b,
\end{align*}
where $q=\lceil4b\rceil$. Since $4b\leq q\leq5b$ for $b\geq1$, taking $TL\sqrt d$ sufficiently small makes the last bound at most $e^{-3/4}$. Raising to the $q$-th power proves $\mu\bigl(|\ell_T|>b\bigr)\leq e^{-3q/4}\leq e^{-3b}$.
The same argument holds at every earlier endpoint.
By Cavalieri's principle, using the trivial probability bound below $b=1$ and the exponential tail above it,
\[
 \mu\bigl(e^{2|\ell_t|}\bigr)
\ =\ 1+2\int_0^\infty e^{2b}\mu\bigl(|\ell_t|>b\bigr)\,db\ \lesssim\ 1.
\]
Since $(\rho_t/\rho)^2\leq e^{2|\ell_t|}$, this proves \eqref{eq:acceptance-heat}.
\end{proof}

To control the conditional velocity fluctuations of the midpoint density, we record Gaussian $L^q$ Poincar\'e in the form
\begin{equation}\label{eq:readable-gaussian-poincare}
 \bigl\|f-\varphi f\bigr\|_{L^q(\varphi)}
\ \lesssim\ \sqrt q\,\bigl\||\nabla f|\bigr\|_{L^q(\varphi)}\quad\text{for all $q\geq2$.}
\end{equation}
For completeness, the Gaussian rotation argument of \cite{Pisier86} proves this directly.
For the coordinate Gaussians $G,G'$ under $\varphi\otimes\varphi$, put $G_\theta=G\cos\theta+G'\sin\theta$ and $G'_\theta=-G\sin\theta+G'\cos\theta$.
Jensen and the fundamental theorem of calculus imply
\begin{align*}
 \bigl\|f-\varphi f\bigr\|_{L^q(\varphi)}\ &\leq\ \bigl\|f(G)-f(G')\bigr\|_{L^q(\varphi\otimes\varphi)}\\
\ &\leq\ \int_0^{\pi/2}\bigl\|\nabla f(G_\theta)\cdot G'_\theta\bigr\|_{L^q(\varphi\otimes\varphi)}\,d\theta\\
\ &\lesssim\ \sqrt q\,\bigl\||\nabla f|\bigr\|_{L^q(\varphi)}.
\end{align*}
The last inequality conditions on $G_\theta$, since $G'_\theta$ is an independent standard Gaussian.
The estimate extends to locally Lipschitz functions by smooth approximation with cutoffs. All functions and weak derivatives used below have polynomial growth.
Applying this componentwise at order two yields
\begin{equation}\label{eq:readable-vector-poincare}
 \mathbb E\bigl|V(G)-\mathbb EV(G)\bigr|^2\ \lesssim\ \mathbb E\bigl\|DV(G)\bigr\|_F^2.
\end{equation}
Also, for $G\sim\varphi$, $\bigl\||G|\bigr\|_{L^q(\varphi)}\lesssim \sqrt{d+q}$, since Gaussian integration by parts evaluates $\mathbb E|G|^{2k}=\prod_{j=0}^{k-1}(d+2j)\leq(d+2k)^k$.

We now express the unnormalized measure $e^{sX}\hat\mu$ in coordinates where the Gaussian remainder estimate will apply. Set $a=h/2$, so $s=a^2/2$, and introduce midpoint coordinates
\[
 w\ =\ v-a\nabla U(x),\qquad y\ =\ x+aw,
 \qquad x\ =\ y-aw,\quad v\ =\ w+a\nabla U(y-aw).
\]
Both shears preserve volume.
For
\[
 B(y,w)\ =\ U(y)-U(y-aw)-aw\cdot\nabla U(y-aw),
 \qquad \mathbb E_0\ =\ \mathbb E_{\mu(dy)\varphi(dw)},
\]
expanding $\bigl|w+a\nabla U(x)\bigr|^2/2$ proves the exact identity
\begin{equation}\label{eq:readable-midpoint-density}
 e^{sX(x)}\mu(dx)\varphi(dv)
\ =\ e^{B(y,w)-s\mu F}\mu(dy)\varphi(dw).
\end{equation}
By \Cref{lem:readable-gibbs} and Jensen, $1\leq Z_s:=\mu\bigl(e^{sX}\bigr)\lesssim1$ when $sL\sqrt d$ is sufficiently small.
Thus $\mathbb E_0e^B=e^{s\mu F}Z_s$, and $\varrho=e^B/\mathbb E_0e^B$ is the density of the image of $Z_s^{-1}e^{sX}\hat\mu$ relative to $\mu\otimes\varphi$.

\begin{lemma}[Midpoint Density]\label{lem:readable-midpoint}
Suppose $\mu$ satisfies the assumptions of \Cref{lem:readable-gibbs}, and let $a\geq0$.
If $a^2L\sqrt d$ is sufficiently small, then, for every integer $q\geq1$,
\[
 \bigl\|B-\mathbb E_0B\bigr\|_{L^q(\mu\otimes\varphi)}\ \lesssim\ qa^2L\sqrt d,
 \qquad \mathbb E_0\varrho^2\ \lesssim\ 1.
\]
\end{lemma}

\begin{proof}
We estimate $B(y,G)-\mathbb E_GB(y,G)$ conditionally on $y$, and bound the variation of its mean with $y$ using Laplacian moments.
Write $P_rf(y)=\mathbb E f(y-rG)=e^{(r^2/2)\Delta}f(y)$. Since $Y-rG$ has density $\rho_{r^2/2}$ when $Y\sim\mu$ and $G\sim\varphi$ are independent, Jensen and H\"older imply
\[
 \|P_rf\|_{L^q(\mu)}^q
\ \leq\ \int |f|^q\frac{\rho_{r^2/2}}\rho\,d\mu
\ \leq\ \left\|\frac{\rho_{r^2/2}}\rho\right\|_{L^2(\mu)}\|f\|_{L^{2q}(\mu)}^q.
\]
Apply this to $f=\Delta U-\mu(\Delta U)$. \Cref{lem:readable-heat}, used up to time $a^2/2$, and \eqref{eq:acceptance-trace} imply
\begin{equation}\label{eq:readable-shifted-trace}
 \bigl\|P_r\Delta U-\mu(\Delta U)\bigr\|_{L^q(\mu)}\ \lesssim\ qL\sqrt d
 \quad\text{for all $0\leq r\leq a$.}
\end{equation}
Gaussian integration by parts relates the gradient and Laplacian averages:
\[
 \mathbb E\bigl[G\cdot\nabla U(y-aG)\bigr]\ =\ -aP_a\Delta U(y),
 \qquad\frac d{dr}P_rU(y)\ =\ rP_r\Delta U(y).
\]
The derivative is dominated on bounded $r$-intervals by a constant depending on $y$ times $1+|G|^2$.
Integrating yields
\begin{equation}\label{eq:readable-B-mean}
 b(y):\ =\ \mathbb E_G B(y,G)
\ =\ a^2P_a\Delta U(y)-\int_0^a rP_r\Delta U(y)\,dr.
\end{equation}
Subtracting $(a^2/2)\mu(\Delta U)$ and applying Minkowski's inequality and \eqref{eq:readable-shifted-trace} proves
\[
 \bigl\|b-(a^2/2)\mu(\Delta U)\bigr\|_{L^q(\mu)}\ \lesssim\ qa^2L\sqrt d,
 \qquad\bigl\|b-\mu b\bigr\|_{L^q(\mu)}\ \lesssim\ qa^2L\sqrt d.
\]
On the other hand, $\nabla_wB=a^2\nabla^2U(y-aw)w$, so \eqref{eq:readable-gaussian-poincare} implies
\[
 \sup_y\bigl\|B(y,G)-b(y)\bigr\|_{L^q(\varphi)}
\ \lesssim\ a^2L\sqrt{q(d+q)}\ \lesssim\ qa^2L\sqrt d.
\]
Order one follows from order two.
The decomposition $B-\mathbb E_0B=(B-b)+(b-\mu b)$ and the triangle inequality prove the moment estimate.
All integrals are finite since $0\leq B(y,w)\leq La^2|w|^2/2$.
Finally, for $B_0=B-\mathbb E_0B$, the moment series, $n^n/n!\leq e^n$ and Jensen's bound $\mathbb E_0e^{B_0}\geq1$ imply
\[
 \mathbb E_0e^{2|B_0|}\ =\ 1+\sum_{n\geq1}\frac{2^n\mathbb E_0|B_0|^n}{n!}\ \lesssim\ 1,
 \qquad
 \mathbb E_0\varrho^2
\ =\ \frac{\mathbb E_0e^{2B_0}}{\bigl(\mathbb E_0e^{B_0}\bigr)^2}\ \lesssim\ 1,
\]
provided $a^2L\sqrt d$ is sufficiently small.
\end{proof}

\subsection{The Corrected Leapfrog Energy Error}
The density estimate from the preceding subsection reduces the weighted remainder moment to a Gaussian moment, uniformly in the midpoint position. We prove that moment using the Frobenius bound, then combine the estimates to complete \Cref{lem:acceptance-local}.
The exact correction \cite[Eq.~(20)]{Bussi07}, equivalently \cite[Lem.~2.1]{Nawaf26persistence}, identifies $\mathcal R$ in midpoint coordinates with
\[
 \mathcal R_{\mathrm{mid}}(y,w)\ =\ U(y+aw)-U(y-aw)
          -aw\cdot\bigl[\nabla U(y+aw)+\nabla U(y-aw)\bigr].
\]
For the corrected remainder, the following Gaussian estimate replaces the third-derivative estimates used to bound the full leapfrog energy error in \cite[Lem.~8 and~12--13, App.~B.2]{chen2023}:
\begin{equation}\label{eq:readable-midpoint-moment}
 \sup_y\bigl\|\mathcal R_{\mathrm{mid}}(y,G)\bigr\|_{L^q(\varphi)}
\ \lesssim\ a^3L_Hq\sqrt{d+q}\quad\text{for all $q\geq2$.}
\end{equation}
For fixed $y$, the remainder is odd in $w$ and
\[
 \nabla_w\mathcal R_{\mathrm{mid}}\ =\ -a^2V(w),
 \qquad V(w)\ =\ \bigl[\nabla^2U(y+aw)-\nabla^2U(y-aw)\bigr]w.
\]
The weak third derivatives satisfy $\bigl\|D\nabla^2U(x)[u]\bigr\|_F\leq L_H|u|$ almost everywhere by the Frobenius Hessian Lipschitz assumption. Their symmetry means that $D\nabla^2U(x)[u]w=D\nabla^2U(x)[w]u$:
\[
 \begin{gathered} DV(w)\ =\ \nabla^2U(y+aw)-\nabla^2U(y-aw)+aD\nabla^2U(y+aw)[w]
                    +aD\nabla^2U(y-aw)[w],\\
 \bigl\|DV(w)\bigr\|_F\ \lesssim\ aL_H|w|. \end{gathered}
\]
Gaussian integration by parts computes the mean, which is bounded using $|\nabla\Delta U|\leq\sqrt d L_H$:
\[
 \mathbb EV(G)\ =\ a\mathbb E\bigl[\nabla\Delta U(y+aG)+\nabla\Delta U(y-aG)\bigr],
 \qquad \bigl|\mathbb EV(G)\bigr|\ \lesssim\ aL_H\sqrt d.
\]
The weak derivatives have polynomial growth, so these integrations and Gaussian Poincar\'e are valid under the assumed $C^2$ regularity.
By \eqref{eq:readable-vector-poincare}, the centered second moment of $V(G)$ is $\lesssim a^2L_H^2d$. Together with the bound on its mean, this implies $\mathbb E\bigl|V(G)\bigr|\lesssim aL_H\sqrt d$.
Since the weak gradient of $|V|$ is bounded by $\|DV\|_F$, \eqref{eq:readable-gaussian-poincare} yields
\[
 \bigl\||V|\bigr\|_{L^q(\varphi)}\ \lesssim\ \mathbb E|V|+\sqrt q\,\bigl\|\|DV\|_F\bigr\|_{L^q(\varphi)}
\ \lesssim\ aL_H\sqrt{q(d+q)}.
\]
A second application to the mean-zero $\mathcal R_{\mathrm{mid}}$ proves \eqref{eq:readable-midpoint-moment}.

\begin{proof}[Proof of \Cref{lem:acceptance-local}]
The endpoint estimate follows from the centered-gradient moments. For the remainder, we use the exact midpoint change of variables and H\"older to combine the density bound with the Gaussian estimate.
By \Cref{lem:readable-gibbs}, $\mu\bigl(e^{2s|X|}\bigr)\lesssim1$ and $1\leq Z_s\lesssim1$.
The inequality $\bigl|e^u-1\bigr|\leq|u|e^{|u|}$ and H\"older prove the stronger first bound
\begin{align*}
 \bigl\|e^{sX/p}-1\bigr\|_{L^p(\hat\mu)}
\ &\leq\ \frac{s}{p}\|X\|_{L^{2p}(\mu)}
                      \mu\bigl(e^{2s|X|}\bigr)^{1/(2p)}\\
\ &\lesssim\ h^2L\bigl(\sqrt{d/p}+1\bigr)\ \lesssim\ h^2L\sqrt{d+p}.
\end{align*}
By the midpoint change of variables \eqref{eq:readable-midpoint-density},
\begin{align*}
 \left(\int e^{sX}|\mathcal R|^p\,d\hat\mu\right)^{1/p}
\ &=\ \left(Z_s\mathbb E_0\bigl[\varrho|\mathcal R_{\mathrm{mid}}|^p\bigr]\right)^{1/p}\\
\ &\leq\ \bigl(Z_s\|\varrho\|_{L^2(\mu\otimes\varphi)}\bigr)^{1/p}\|\mathcal R_{\mathrm{mid}}\|_{L^{2p}(\mu\otimes\varphi)}
\\\ &\lesssim\ h^3L_Hp\sqrt{d+p}.
\end{align*}
Here \Cref{lem:readable-midpoint} bounds the density in $L^2(\mu\otimes\varphi)$, and \eqref{eq:readable-midpoint-moment} bounds the Gaussian moment uniformly in $y$.
This proves \Cref{lem:acceptance-local}.
All smallness conditions follow by taking $h^2L\sqrt d$ sufficiently small.
\end{proof}

\end{document}